\documentclass[10pt]{article}

\usepackage[includeheadfoot,top=2 cm, bottom=2 cm, left=2 cm, right=2 cm]{geometry}
\usepackage{graphicx}
\usepackage{lipsum}
\usepackage{amsfonts,amsthm}
\usepackage{amssymb}
\usepackage{upgreek}
\usepackage{hyperref}
\usepackage{algorithmicx,algorithm,algcompatible,eqparbox}

\usepackage{multicol,tabularx,capt-of}
\makeatletter
\renewcommand{\ALG@beginalgorithmic}{\small}
\makeatletter
\usepackage{multirow}
\usepackage{array}
\usepackage{mathtools}
\usepackage{color}
\usepackage{xcolor}
\usepackage{float}
\usepackage{subcaption}
\usepackage[noabbrev,capitalise]{cleveref}
\hypersetup{
	colorlinks,
	linkcolor={blue!},
	citecolor={blue!},
	urlcolor={blue!}
}

\usepackage[
sortcites,
backend=bibtex,
hyperref=true,
firstinits=true,
maxbibnames=200,
backref,
backrefstyle=none
]{biblatex} 
\renewbibmacro{in:}{}
\Crefname{proposition}{Proposition}{Propositions}
\crefname{equation}{}{}
\Crefname{equation}{}{}

\newcommand{\cond}{\textup{cond}}

\newcommand{{{\bxi}}}{\mathbf{\upxi}}
\newcommand{\Frob}{\mathrm{F}}

\newcommand{\ba}{\mathbf{a}}

\newcommand{\bb}{\mathbf{b}}

\newcommand{\bx}{\mathbf{x}}

\newcommand{\bTheta}{\mathbf{\Theta}}
\newcommand{\bSigma}{\mathbf{\Sigma}}

\newcommand{\bPhi}{\mathbf{\Phi}}
\newcommand{\bR}{\mathbf{R}}
\newcommand{\bA}{\mathbf{A}}
\newcommand{\bE}{\mathbf{E}}

\newcommand{\bB}{\mathbf{B}}
\newcommand{\bhS}{\widehat{\mathbf{S}}}
\newcommand{\bhR}{\widehat{\mathbf{R}}}

\newcommand{\bhQ}{\widehat{\mathbf{Q}}}

\newcommand{\bhZ}{\widehat{\mathbf{Z}}}

\newcommand{\bhP}{\widehat{\mathbf{P}}}
\newcommand{\bhY}{\widehat{\mathbf{Y}}}

\newcommand{\bX}{\mathbf{X}}
\newcommand{\bZ}{\mathbf{Z}}

\newcommand{\bD}{\mathbf{D}}

\newcommand{\bW}{\mathbf{W}}

\newcommand{\bGamma}{\mathbf{\Gamma}}

\newcommand{\bS}{\mathbf{S}}
\newcommand{\bY}{\mathbf{Y}}
\newcommand{\bQ}{\mathbf{Q}}

\newcommand{\bP}{\mathbf{P}}
\newcommand{\bI}{\mathbf{I}}

\newcommand{\bU}{\mathbf{U}}
\newcommand{\bL}{\mathbf{L}}
\newcommand{\bV}{\mathbf{V}}

\newcommand{\bPi}{\mathbf{\Pi}}

\newtheorem{lemma}{lemma}[section]
\newtheorem{proposition}[lemma]{Proposition}
\newtheorem{corollary}[lemma]{Corollary}
\newtheorem{theorem}[lemma]{Theorem}
\newtheorem{remark}[lemma]{Remark}
\newtheorem{definition}[lemma]{Definition}
\newtheorem{assumptions}[lemma]{Assumptions}

\begin{document}
	\title{Randomized Cholesky QR factorizations}
	\author{Oleg Balabanov\thanks{Part of this work was conducted while the author was at Sorbonne Universit\'e, Inria, CNRS, Universit\'e de Paris, Laboratoire Jacques-Louis Lions, Paris, France. Email: oleg.balabanov@inria.fr.}}		
	\date{}
	\maketitle
	
	\begin{abstract}
		This article proposes and analyzes several variants of the randomized Cholesky QR factorization of a matrix $X$. Instead of computing the R factor from $X^T X$, as is done by standard methods, we obtain it from a small, efficiently computable random sketch of $X$, thus saving computational cost and improving numerical stability. The proposed direct variant of the randomized Cholesky QR requires only half the flops and the same communication cost as the classical Cholesky QR. At the same time, it is more robust since it is guaranteed to be stable whenever the input matrix is numerically full-rank.  The rank-revealing randomized Cholesky QR variant has the ability to sort out the linearly dependent columns of $X$, which allows to have an unconditional numerical stability and reduce the computational cost when $X$ is rank-deficient.  
		We also depict a column-oriented randomized Cholesky QR that establishes the connection with the randomized Gram-Schmidt process, and a reduced variant that outputs a low-dimensional projection of the Q factor rather than the full factor and therefore yields drastic computational savings.
		It is shown that performing minor operations in higher precision in the proposed algorithms can allow stability with working unit roundoff independent of the dominant matrix dimension. This feature may be of particular interest for a QR factorization of tall-and-skinny matrices on low-precision architectures.
	\end{abstract}
	
	\begin{keywords}
		QR factorization, randomization, sketching, Cholesky, rank-revealing, numerical stability, rounding errors, loss of orthogonality, multi-precision arithmetic, communication-avoiding algorithms.
	\end{keywords}
	
	\section{Introduction} 
	This work is devoted to computing a thin QR factorization 
	$ \bX = \bQ \bR$
	of matrix $\bX \in \mathbb{R}^{m \times n}$ with $n \ll m$, where $\bR$ is upper triangular or trapezoidal, possibly with permuted columns, and $\bQ$ is approximately orthonormal or very well-conditioned.
	Such QR factorizations constitute the basic kernels for many scientific computing algorithms. Applications include solution of least-squares problems~\cite{golub2013matrix}, computation of low-rank approximations~\cite{halko2011algorithm}, solution of linear systems and eigenvalue problems~\cite{golub2013matrix,saad2011numerical,saad2003iterative}, model order reduction~\cite{haasdonk2008reduced}, and more. There is a vast variety of algorithms for computing a QR factorization. They can be majorly divided into three main categories: the ones based on Householder transformations or Givens rotations~\cite{golub2013matrix,higham2002accuracy}, the ones based on Gram-Schmidt orthogonalization~\cite{leon2013gram}, and the ones based on Cholesky QR~\cite{fukaya2014choleskyqr2}. We are here interested in the latter type of algorithms.  
	
	In scientific computing, special attention must be paid to the numerical stability of algorithms, that is, their sensitivity to rounding errors in finite precision arithmetic~\cite{higham2002accuracy}. A computed QR factorization of $\bX$ can be said to be numerically stable if $\mathrm{cond}(\bQ) = \mathcal{O}(1)$ and the columns of $\bQ \bR$ approximate the columns of $\bX$ up to machine precision.

	\subsection{Cholesky QR} \label{standCholeskyQR}
	A Cholesky QR (CholeskyQR) factorization of $\bX$ proceeds by first obtaining the R factor through a Cholesky factorization of the Gramian $\bX^\mathrm{T} \bX$, and then retrieving the Q factor by forward substitution, as shown in~\cref{alg:CholeskyQR}.
	\begin{algorithm}[h] \caption{Cholesky QR (CholeskyQR)} \label{alg:CholeskyQR}
		\begin{algorithmic}
			\STATE{\textbf{Input:}}
			\STATE{\;\;\; \makebox[0.5cm]{$\bX$} is $m \times n$ matrix} 
			\STATE{\textbf{Output}:}
			\STATE{\;\;\; \makebox[0.5cm]{$\bQ$} is $m \times n$ orthonormal Q factor}
			\STATE{\;\;\; \makebox[0.5cm]{$\bR$} is $n \times n$ upper triangular R factor}
			\STATE{\textbf{function} $[\bQ, \bR]=\mathtt{CholeskyQR}(\bX)$}
			\STATE{1. $\bA \leftarrow \bX^\mathrm{T} \bX$}
			\STATE{2. $\bR \leftarrow \mathtt{chol}(\bA)$}
			\STATE{3. $\bQ \leftarrow  \bQ \bR^{-1}$}
		\end{algorithmic}
	\end{algorithm}
	This procedure is well suited for distributed computing because it requires only one global synchronization between processors. It can be seen as an alternative to another communication-avoiding algorithm called TSQR~\cite{demmel2012communication}. The computational advantages of CholeskyQR over TSQR include the twice less computational cost in terms of flops, lower amount of global synchronizations between processors, and a simpler reduction operator. However, unlike TSQR, which is ultimately stable, CholeskyQR often introduces instabilities that limit its use. In particular, one can guarantee numerical stability of the algorithm only by making sure that the unit roundoff satisfies $u \leq \mathrm{cond}(\bX)^{-2} F(m,n)^{-1}$, where $F(m,n)$ is a low-degree polynomial~\cite{yamamoto2015roundoff}. To alleviate this drawback, a shifted Cholesky QR with reorthogonalization (shifted CholeskyQR2) has recently been proposed, which is stable under the condition that $u \leq \mathrm{cond}(\bX)^{-1}  F(m,n)^{-1}$~\cite{fukaya2020shifted}. The shifted CholeskyQR2, however, doubles the cost of the classical CholeskyQR. Moreover, the aforementioned stability guarantees may still be insufficient for badly scaled or ill-conditioned problems.
	
	This article proposes randomized Cholesky QR factorizations that can be more than four times as efficient as shifted CholeskyQR2 in terms of flops, and twice as efficient in terms of communication cost, and yet they achieve at least as good stability as that of shifted CholeskyQR2 or even unconditional stability like that of Householder QR or TSQR.

	\subsection{Randomized Cholesky QR} 
	We provide two methods to compute a QR factorization: the (direct) randomized Cholesky QR (RCholeskyQR) and the rank-revealing randomized Cholesky QR (RRRCholeskyQR) that are based on a dimension reduction technique called random sketching. See~\cite{martinsson2020randomized,woodruff2014sketching} for an overview of this technique for scientific computing. Both RCholeskyQR and RRRCholeskyQR first compute the factor $\bR$ by a low-dimensional QR factorization of a random sketch $\bTheta \bX \in \mathbb{R}^{k \times n}$ of $\bX$ and then retrieve the factor $\bQ$ by forward substitution. The matrix $\bTheta \in \mathbb{R}^{k \times m}$ is a suitable random matrix typically with $k = \mathcal{O}(n)$ (say $k=2n$) rows that can be efficiently applied to $\bX$ in the given architecture, and that is with high probability an approximate isometry for the column space of $\bX$.  
	In the algorithms, the R factor can be computed for instance with the Cholesky factorization of the sketched Gramian $(\bTheta \bX)^\mathrm{T} (\bTheta \bX)$. In this case the proposed QR factorizations become aligned with the direct and rank-revealing generalized Cholesky QR associated with the sketched inner product $\langle \bTheta \cdot, \bTheta \cdot \rangle_2$. This fact  gives our methods their names. However, since $\bTheta\bX$ is small, computing the R factor from $\bTheta\bX$ using the CholeskyQR may become inappropriate and should be done by more expensive but more reliable methods based, for instance, on Householder transformations or Givens rotations. See~\cref{randCholeskyQR,rrrandCholeskyQR} for more details. It has to be underlined that although randomization entails a possible failure of the algorithms, here this should not cause any concerns. We take the probability of failure as a user-defined parameter that can be chosen to be astronomically small, say $10^{-10}$ or $10^{-20}$ without much impact on computational cost.

	We show that the RCholeskyQR is stable when the matrix $\bX$ is numerically full-rank, which makes it at least as robust as shifted CholeskyQR2. Yet it should take significantly less flops and communication between distributed processors.
	Moreover, the gain in the computational cost can be drastic when only a small projection of $\bQ$ is needed rather than the full matrix (see~\cref{rRColQR}). The RRRCholeskyQR variant is shown to be even more stable than RCholeskyQR, and provides QR factorization with well-conditioned Q factor and column-wise approximation error close to machine precision under the condition $u <  F(m,n)^{-1}$ that does not involve $\bX$. At the same time, this factorization can inherit the high efficiency of RCholeskyQR or even improve it by a great amount if (normalized) $\bX$ is a low-rank matrix. This makes RRRCholeskyQR a very desirable alternative to Householder QR, TSQR and other unconditionally stable QR factorizations in applications.
	It is then noticed that by performing the minor operations in RCholeskyQR and RRRCholeskyQR algorithms in higher precision, one can guarantee stability of the algorithms using working unit roundoff $u$ independent of the high dimension $m$ i.e. to have $F(m,n) = F(n)$. This can be particularly useful for the QR factorization of tall-and-skinny matrices on low precision architectures. 
	
	The RCholeskyQR and RRRCholeskyQR methods are closely related to the Randomized Gram-Schmidt (RGS) QR factorization and its block version recently proposed in~\cite{balabanov2021randomizedGS,balabanov2021randomized}. The RGS process was successfully applied to improving Krylov methods for the solution of linear systems and eigenvalue problems. See~\cite{MATLABbalabanov2022} for an open-source implementation of the RGS-based Krylov methods that can provide more than $10\times$ speedups over the standard built-in functions. Since RGS also proceeds with the orthogonalization of the random sketch of $\bX$, it inherits some of the properties of randomized Cholesky QR algorithms from this paper.

	\paragraph{Related work.} The idea of using a QR factorization with respect to the sketched inner product goes back to~\cite{balabanov2019randomized,balabanov2019randomized2}, where it was employed to improve conditioning of reduced basis for (parametric) linear systems. This approach is aligned with the reduced RCholeskyQR given in~\cref{rRColQR}. In this context, our main contribution is to obtain rigorous guarantees of numerical stability.

	A QR factorization with the RGS algorithm was introduced in~\cite{balabanov2021randomizedGS} and was recently extended to block matrices in~\cite{balabanov2021randomized}. In exact arithmetic, the output of RCholeskyQR and RGS is the same. These algorithms are also closely related from a numerical point of view. As pointed out in~\cite[Remark 2.10]{balabanov2021randomizedGS}, in RGS it is possible to compute the sketch of the newly obtained column of $\bQ$ from the sketches of previously computed vectors, rather than by multiplying that column by $\bTheta$.
	This modification leads to an algorithm that is nothing more than a column-oriented variant of RCholeskyQR (see~\cref{colRCholeskyQR}). However, as noted in~\cite[Remark 2.10]{balabanov2021randomizedGS} and also verified in our numerical experiments, the column-oriented RCholeskyQR may be less numerically stable than RGS. A similar picture is observed for the blockwise version of RCholeskyQR~\cite{balabanov2021randomized}. Nevertheless, we show that RCholeskyQR has stability guarantees similar to RGS when $\bX$ is numerically full-rank, which is a fairly common situation. More information about the connection between the column-oriented RCholeskyQR and RGS can be found in~\cref{colRCholeskyQR}.
	In addition, in~\cite[Section 2.2.2]{balabanov2021randomized} the authors proposed a direct RCholeskyQR algorithm, the same as in this article. However, they omitted its analysis.

	There is another recent work~\cite{fan2021novel}, which in parallel to~\cite{balabanov2021randomized} explored the ideas of the direct RCholeskyQR. In particular, \cite[Algorithm 2.6]{fan2021novel} is similar to the RCholeskyQR2 algorithm presented in this paper. The key difference is that in~\cite{fan2021novel} the authors mainly rely on uniform sampling of rows of $\bX$ rather than sketching with oblivious subspace embeddings or nonuniform sampling. Note that sketching using Gaussian matrices and nonuniform sampling have also been mentioned, but only as a theoretical supplement to the uniform sampling approach. In particular, Gaussian matrices have been left out of the scope of numerical experiments, theoretical discussions of computational costs, and stability analysis. At the same time, the provided in~\cite{fan2021novel} nonuniform sampling is infeasible for practical use, as it requires a priori knowledge of the Q factor. Our article in this sense is more complete. It fully explores sketching with Gaussian matrices both theoretically and numerically. In addition, it takes into account other dimension reduction maps. For instance, we consider sketching with subsampled randomized Hadamard transform, which, unlike the methodology from~\cite{fan2021novel}, leads to provably stable and accurate algorithms that have less complexity and memory consumption compared to shifted CholeskyQR3 or other deterministic algorithms. A similar result can be achieved using nonuniform row sampling approaches from~\cite{drineas2012fast,woodruff2014sketching,clarkson2017low}, which are also considered.
	
	The proposed RRRCholeskyQR method is connected to the randomized QR with column pivoting or the randomized rank-revealing QR from~\cite{martinsson2015blocked,martinsson2017householder,xiao2017fast,duersch2017randomized,mary2015performance}. First of all, here we primarily focus on the QR factorization of tall-and-skinny matrices, and not on the QR factorization of large (low-rank) matrices, considered in the above-mentioned works. Though, the RRRCholeskyQR factorization in principle can also be used in the latter context. 
	The main advantage of RRRCholeskyQR in this case is that, unlike other methods, it calculates the R factor solely from the random sketch, which can significantly reduce computational cost. In particular, RRRCholeskyQR can require up to $\frac{n}{r}$ less flops, where $r$ is the approximation rank.
	
	The idea of performing a CholeskyQR with multi-precision arithmetic was explored in~\cite{yamazaki2015mixed}. In contrast to~\cite{yamazaki2015mixed}, here we propose to increase the rounding accuracy only for minor operations with little impact on the overall computational cost. Moreover, we prove numerical stability for working rounding independent of the dominant dimension $m$. As far as we know, this property is inherent only in randomized algorithms.

	The rest of the paper is organized as follows. \Cref{prem} describes main notations. 
	In~\cref{RS}, we introduce the random sketching technique along with some of the results underlying our randomized algorithms.  In~\cref{randCholeskyQR} we introduce the direct RCholeskyQR factorization and discuss its derivatives, which are a column-oriented version and a reduced version. \Cref{rrrandCholeskyQR} is devoted to RRRCholeskyQR factorization. Numerical stability of the proposed algorithms is characterized in~\cref{stability}. The methodology is validated numerically in~\Cref{numexperiments}. \Cref{conclusion} concludes the article.

	\subsection{Preliminaries} \label{prem}
	In the manuscript, we use the following notations and assumptions. For a given matrix $\bA$, the $i$-th column of $\bA$ is denoted by $\bA_{(i)}$ or $\bA_{(:,i)}$, and the $i$-th row is denoted by $\bA_{(i,:)}$. A submatrix consisting of the $i$-th through the $j$-th consecutive columns of $\bA$ is denoted by $\bA_{(i:j)}$ or $\bA_{(:,i:j)}$, and the one that consists of the $i$-th through the $j$-th consecutive rows, by $\bA_{(i:j,:)}$. In addition, the $k$-th column of the submatrix $\bA_{(i:j,:)}$ is denoted by $\bA_{(i:j,k)}$. If $\bA$ is a block matrix, then the above notation is used to denote submatrices composed of the corresponding blocks of $\bA$. For instance, in this case $\bA_{(i:j,k)}$ denotes a submatrix composed by column-wise concatenation of $(i,k)$-th, $(i+1,k)$-th, $\hdots$, $(j,k)$-th blocks of $\bA$. We denote by $|\bA|$ the matrix whose entries are the absolute values of the corresponding entries of the matrix $\bA$. We let $\sigma_{min}(\bA)$ and $\sigma_{max}(\bA)$ denote the minimal and maximal singular values of $\bA$, and $\mathrm{cond}(\bA)$ to denote the condition number $\frac{\sigma_{max}(\bA)}{\sigma_{min}(\bA)}$. In addition, $\bA^\mathrm{T}$ and $\bA^\mathrm{\dagger}$ denote the transpose of $\bA$ and the Moore-Penrose pseudo-inverse of $\bA$, respectively. We let $\bI$ be an identity matrix of size appropriate for the expression in which it is used. An arithmetic expression or quantity $A$ computed using finite precision arithmetic is denoted by $\widehat{A}$ or $fl(A)$. In addition, in the article we will assume that the dominant operations in randomized algorithms are performed with unit roundoff $u$, and the secondary operations, which are random projections and low-dimensional operations, use unit roundoff $u_{f}$ lower than $u$ by sufficiently large low-degree polynomial in $m$ and $n$.

	\section{Random sketching technique} \label{RS}
	
	Let $\bTheta \in \mathbb{R}^{k \times m}$ be a sketching matrix with $k \ll m$ rows.  This matrix is seen as an embedding of low-dimensional subspaces of $\mathbb{R}^m$ into $\mathbb{R}^{k}$. In addition, it is chosen such that it is an approximate isometry of the subspace(s) of interest, or in other words, an $\varepsilon$-embedding.
	
	\begin{definition} \label{def:subemb}
		Sketching matrix $\bTheta$ is called an $\varepsilon$-embedding for $\bV \in \mathbb{R}^{m \times d}$ or $\mathrm{range}(\bV)$ if 
		$$ (1-\varepsilon)\|\bV \bx \|^2_2 \leq \| \bTheta \bV \bx \|^2_2  \leq (1+\varepsilon) \|\bV \bx \|^2_2.$$
	\end{definition}
	Our analysis assumes that the parameter $\varepsilon$ is set by the user. However, it should be underlined that for the applications in this paper it is not necessary to consider very small values of $\varepsilon$. It suffices to take $\varepsilon = \frac{1}{2}$ or $\varepsilon = \frac{1}{4}$.
	
	As pointed out in~\cite{balabanov2019randomized,balabanov2021randomizedGS}, the $\varepsilon$-embedding property of $\bTheta$ allows $\bV$ to be approximately orthonormalized by orthonormalizing the sketch $\bTheta \bV$. This observation follows from~\cref{thm:epscond} and underlies all our randomized algorithms. 
	\begin{corollary}[Corollary 2.2 in~\cite{balabanov2021randomizedGS}] \label{thm:epscond}
		Let $\bV \in \mathbb{R}^{m \times d}$ be some matrix.
		If $\bTheta$ is {an} $\varepsilon$-embedding for $\bV$, then the singular values of $\bV$ are bounded by $$ (1+\varepsilon)^{-\frac{1}{2}} \sigma_{min}(\bTheta\bV) \ \leq \sigma_{min}(\bV) \leq \sigma_{max}(\bV) \leq  (1-\varepsilon)^{-\frac{1}{2}} \sigma_{max}(\bTheta\bV).$$
	\end{corollary}
	
	In algorithms, we prefer to construct $\bTheta$ without any a priori information about the $\bV$ matrix. This can be done in a probabilistic way by drawing $\bTheta$ from a carefully designed probability distribution such that $\bTheta$ satisfies the $\varepsilon$-embedding property for any fixed $m \times d$ matrix with high probability. Such $\bTheta$ will be called an ($\varepsilon$,$\delta$,$d$)-oblivious subspace embedding (or ($\varepsilon$,$\delta$,$d$)-OSE), as defined below.
	\begin{definition} \label{def:OSE}
		Random matrix $\bTheta$ is called an ($\varepsilon$,$\delta$,$d$)-OSE if it is an $\varepsilon$-embedding for any (fixed) $m \times d$ matrix $\bW$ with probability at least $1-\delta$.
	\end{definition}
	
	The advantage of our randomized algorithms is that they do not rely on a particular OSE distribution for $\bTheta$ but rather allow the ability to choose an OSE depending on the given application and the computational architecture to obtain the most computational benefit. For instance, in a classical sequential environment, the most beneficial OSE can be the subsampled randomized Hadamard transform (SRHT) as it can be applied to vectors with as little as $n \log_2 n $ flops (compared to the $2nk$ required by unstructured matrices) due to the tensor structure of the Hadamard transform. On the other hand, if memory consumption is the main concern, then both structured and unstructured OSEs can result in significant computational savings since they can be constructed and operated with a seeded random number generator, implying negligible memory consumption of $\bTheta$. The advantage of (rescaled) Gaussian or Rademacher OSEs over SRHT is their high suitability for cache-based or massively parallel computational environments. Finally, the CountSketch OSE is computationally advantageous if the $\bV$ matrix we want to orthogonalize is sparse. In this work, we chose SRHT and rescaled Gaussian matrices as representative OSEs. The SRHT matrix is defined as a product of a diagonal matrix of random signs (possibly padded with zeros to make the output dimension a power of $2$) with structured Hadamard matrix, followed by a uniform sampling matrix and a scaling factor $\textstyle \sqrt{\frac{1}{k}}$. It follows from~\cite{balabanov2019randomized,tropp2011improved,boutsidis2013improved} that SRHT matrix with 
	\begin{subequations} \label{eq:OSEsdim}
		\begin{equation}
		k \geq 2( \varepsilon^{2} - \varepsilon^3/3)^{-1} \left (\sqrt{d}+ \sqrt{8 \log{\textstyle \frac{6m}{\delta}}} \right )^2 \log{\textstyle \frac{3d}{\delta}} 
		\end{equation}
		rows is an ($\varepsilon$,$\delta$,$d$)-OSE. The rescaled Gaussian matrices have entries that are i.i.d. Gaussian variables, scaled by a factor $\sqrt{\frac{1}{k}}$. Such matrices satisfy the OSE property, if~\cite{balabanov2019randomized}
		\begin{equation}
		k \geq 7.87 \varepsilon^{-2}( {6.9} d + \log{\textstyle \frac{1}{\delta}}).
		\end{equation}
	\end{subequations}
	We see that for both SRHT and Gaussian matrices the required number of rows $k$ is independent or only logarithmically dependent on the dimension $m$ and the probability of failure $\delta$, suggesting the potential dimension reduction with these OSEs. Furthermore, it has to be noted that even-though the theoretical bounds for SRHT matrices are somewhat worse than for Gaussian matrices, in our applications these embeddings provide practically the same results. In particular, it is revealed that for both SRHT and Gaussian matrices the sampling dimension $k=\mathcal{O}(d)$ (say $k=2d$ or $k=4d$) should be sufficient.  
	
	Apart from OSEs, there is another way to efficiently obtain an $\varepsilon$-embedding that should be mentioned. However, it would require a priori knowledge of the $\bV$ matrix, and therefore would not have all the attendant advantages that the OSEs have. We can proceed as follows. We can choose the sketching matrix $\bTheta$ to define a rescaled non-uniform sampling according to a probability distribution $q_1, q_2, \hdots q_m$, the so-called leverage-scores, as explained in~\cref{levscores}.
	
	\begin{proposition}[Leverage scores sampling, corollary of Theorem 17 in~\cite{woodruff2014sketching}] \label{levscores}
		Let $\bV$ be some $m \times d$ matrix. 
		Let $\bTheta = \bGamma \bD$, where $\bGamma$ is $k \times m$ sampling with replacement matrix that samples the $i$-th entry of the input vector with probability $q_i$, and $\bD$ is $m \times m$ diagonal matrix with entries $\bD_{(i,i)} = \frac{1}{\sqrt{k q_i}}$. Furthermore, assume that $q_1,q_2, \hdots q_m$ are such that    
		$$ q_i \geq \beta \frac{1}{k} \| \bW_{(i,:)} \|^2_2,$$
		where $\beta>0$ is some parameter, and $\bW_{(i,:)}$ are rows of some orthonormal matrix $\bW$ with $\mathrm{range}(\bW) = \mathrm{range}(\bV)$, $1 \leq i \leq m$. 	
		If 
		$$ k >  144 d \varepsilon^{-2} \ln{\textstyle \frac{2d}{\delta}}/\beta,$$
		then $\bTheta$ is an $\varepsilon$-embedding for $\bV$ with probability at least $1-\delta$.
	\end{proposition}
	
	For a given $\bV$, one can calculate the distribution $q_1, q_2, \hdots q_m$ with the approximation parameter $\beta = \textstyle \frac{1}{2}$ in~\cref{levscores} using only $\mathcal{O}(md \log d + d^3)$ flops~\cite{drineas2012fast} or even $\mathcal{O} ( \mathrm{nnz}(\bV) \log d + d^3)$ flops~\cite{clarkson2017low}. Consequently, in this way one can efficiently obtain a sampling matrix $\bTheta$ with $k = \mathcal{O}(d \log d)$ rows, which satisfies the $\varepsilon$-embedding property for $\bV$ with high probability. We omit further details.
	
	Furthermore, it was shown in~\cite{balabanov2021randomizedGS} that sketching with an OSE should not significantly increase rounding errors (in the worst case). This important result will be used in~\cref{stability} to characterize the stability of the proposed randomized algorithms using rounding that does not depend on the dimension $m$. It is summarized in~\cref{thm:thetadeltax2}.
	\begin{corollary}[Corollary of~Theorem 2.5 and Corollary 2.6 in~\cite{balabanov2021randomizedGS}] \label{thm:thetadeltax2}
		Given $\bY \in \mathbb{R}^{m\times n}$ and $\bZ \in \mathbb{R}^{n \times l} $ possibly depending on $\bTheta$, consider the product 
		$$\bY \bZ,$$
		computed with finite precision arithmetic. Consider {probabilistic} rounding model, where the rounding errors due to each elementary arithmetic operation are random variables possibly depending on each other, but are independently centered. Furthermore assume that the errors are bounded so that, it holds 
		$$|\bY \bZ - fl(\bY \bZ)| \leq \bU,$$
		for some matrix $\bU$.
		If $\bTheta$ is a $(\varepsilon/4,l^{-1}\binom{m}{d}^{-1}\delta, d)$ OSE, with $d = 4.2 c^{-1} \log{\textstyle \frac{4}{\delta}}$, where $c \leq 1$ is some universal constant, then 
		\begin{equation}\label{eq:thetadeltax2}
		\| \bTheta(\bY \bZ_{(:,i)} - fl(\bY \bZ_{(:,i)})) \|_2 \leq \sqrt{1+\varepsilon} \|\bU_{(:,i)} \|_2
		\end{equation}
		holds with probability at least $1-2\delta$ for $i = 1,2, \hdots, l$ simultaneously. 
	\end{corollary}
	\Cref{thm:thetadeltax2} says that in practice the sketch $\bTheta(\bY \bZ - fl(\bY \bZ))$ of the rounding error matrix should have column norms not much larger than the worst-case bound of $\bY \bZ - fl(\bY \bZ)$. We notice an improvement by nearly a factor of $\sqrt{\frac{m}{k}}$ over the following trivial estimate (for SRHT matrices):
	
	$$ \| \bTheta(\bY \bZ_{(:,i)} - fl(\bY \bZ_{(:,i)})) \|_2 \leq \|\bTheta \|_2 \|\bY \bZ_{(:,i)} - fl(\bY \bZ_{(:,i)})\|_2 = {\textstyle\sqrt{\frac{m}{k}}} \|\bU_{(:,i)}\|_2,$$
	$i = 1,2 \hdots, l$. The condition of~\cref{thm:thetadeltax2} for $\varepsilon = \frac{1}{2}$ and $l \leq m$ can be satisfied by a Gaussian OSE with $\mathcal{O}(\log{m}\log{{\textstyle \frac{1}{\delta}}})$ rows. For SRHT, this requirement is $\mathcal{O}(\log^2{m}\log^2{\textstyle \frac{1}{\delta}})$, although, as has been said, in practice SRHT and Gaussian matrices give very similar results.
	
	\section{Randomized Cholesky QR} \label{randCholeskyQR}
	The efficiency and stability of CholeskyQR of $\bX$ can be improved by changing the $\ell_2$-orthogonality condition 
	$$ \bQ^\mathrm{T} \bQ = \bI,$$
	to the sketched one
	\begin{equation}~\label{eq:sketchedcond}
	(\bTheta \bQ)^\mathrm{T} (\bTheta \bQ) = \bI,
	\end{equation}
	where $\bTheta$ is an $\varepsilon$-embedding for $\bX$. The sketching matrix $\bTheta \in \mathbb{R}^{k \times m}$ can be readily taken as a low-dimensional OSE or a leverage score sampling matrix described in~\cref{RS}. Furthermore, in our experiments, it is revealed that the required theoretical bounds~\cref{eq:OSEsdim} for OSEs are pessimistic. In our applications, using SRHT or Gaussian matrices with just $k=2n$ rows should be sufficient. 
	
	According to~\cref{thm:epscond}, the matrix $\bQ$ that satisfies~\cref{eq:sketchedcond} is very well-conditioned. The obtained sketched QR factorization can be used directly in randomized methods such as block RGS process from~\cite{balabanov2021randomized} or sketched Galerkin and minres approximations with a reduced basis~\cite{balabanov2019randomized,balabanov2019randomized2}. Alternatively, if having well-conditioned Q factor is insufficient, the sketched QR can be post-processed by the classical CholeskyQR to obtain a Q factor orthonormal to machine precision.
	\subsection{Direct RCholeskyQR}
	\Cref{alg:RCQR} depicts RCholeskyQR algorithm for computing a QR factorization of $\bX$ that satisfies~\cref{eq:sketchedcond}. In step 2 the algorithm computes a QR factorization of a small matrix $\bP$. This task can be performed for instance with the classical CholeskyQR in sufficient precision. In this case~\cref{alg:RCQR} can be viewed as the generalized Cholesky QR with respect to the sketched (randomized) inner product, which gives~\cref{alg:RCQR} its name. However, a better way can be to appeal to more stable methods in step 2 such as Householder QR or Givens QR. Note that since the matrix $\bP$ is small, using more efficient or less efficient method in step 2 should not have much impact on the overall computational cost of~\cref{alg:RCQR}.

	It is shown in~\cref{stabRCholeskyQR} that RCholeskyQR is stable whenever $\bX$ is numerically full-rank i.e. when 
	\begin{equation} \label{eq:condstab}
	\mathrm{cond}(\bX^*) \leq F(m,n)^{-1} u^{-1}, 
	\end{equation}
	where $F(n,m)$ is a low-degree polynomial, $\bX^*$ is $\bX$ with normalized columns, and $u$ is the unit roundoff. Note that this condition implies that RCholeskyQR should be at least as stable as shifted CholeskyQR2, or even more stable if $\bX$ is ill-conditioned or has large first dimension $m$. Our experiments in~\cref{numexperiments} confirm this fact. Furthermore, by performing non-dominant operations in higher precision we can use $F(m,n)$ in~\cref{eq:condstab} independent of $m$. 
	
	%\begin{remark}
	%	RCholeskyQR may fail if~\cref{eq:condstab} does not hold. This problem can be alleviated by slight random perturbation of $\bX$, i.e. by multiplying each entry $\bX_{(i,j)}$ by $1+\tau g_{i,j}$, where $g_{i ,j}$ are Gaussian random variables and $\tau$ is a tolerance parameter close to machine precision. This may improve the conditioning of $\bX$, and yet should only add a small error to the $\bQ\bR$ approximation. It is suggested by~\cite{sankar2006smoothed,azais2004upper} that in this case the condition number of (column-normalized) $\bX$ is bounded with probability at least $1-\mathcal{O}(\frac{u}{\tau})$. Such a guarantee however can be insufficient in applications. We leave this topic for future research.
	%\end{remark}

	Let us now characterize the performance of~\cref{alg:RCQR} {in} different computational architectures. Using the SRHT matrix as $\bTheta$, the computational cost of RCholeskyQR in a classical sequential environment is mainly determined by the computation of $\bX \bR^{-1}$ with forward substitution, which in total requires $mn^2$ flops. 
	On distributed architectures, RCholeskyQR should consume only one global synchronization between processors. Moreover, if $\bTheta$ is an unstructured OSE, the sketching step is an explicit matrix-matrix product that can be performed with level 3 BLAS routine. In this case the reduction operator is a simple addition as in standard CholeskyQR. This suggests that RCholeskyQR should be four times more efficient in terms of flops, and twice as efficient in terms of communication cost as standard/shifted CholeskyQR2. Note also that both Gaussian and SRHT matrices should have a negligible storage cost due to the use of a seeded random number generator.	
	The sketching step $\bP \leftarrow \bTheta \bX$ in RCholeskyQR requires only one pass over $\bX$, just like computing the Gramian $\bX^\mathrm{T} \bX$ in standard CholeskyQR. From this we conclude that RCholeskyQR should have half the cost of a standard/shifted CholeskyQR2 in terms of data passes. In addition, RCholeskyQR can be even more advantageous when $\bX$ is stored column-wise, since in this case the computation of $\bTheta \bX$ can still be done in single pass, while the computation of $\bX^\mathrm{T} \bX$ cannot.
	\begin{algorithm}[h] \caption{Randomized Cholesky QR (RCholeskyQR)} \label{alg:RCQR}
		\begin{algorithmic}
			\STATE{\textbf{Input:}}
			\STATE{\;\;\; \makebox[0.5cm]{$\bX$} is $m \times n$ matrix} 
			\STATE{\;\;\; \makebox[0.5cm]{$\bTheta$} is $k \times m$ sketching matrix (possibly provided as a function handle)} 		
			\STATE{\textbf{Output}:}
			\STATE{\;\;\; \makebox[0.5cm]{$\bQ$} is $m \times n$ well-conditioned Q factor}
			\STATE{\;\;\; \makebox[0.5cm]{$\bS$} is $k \times n$ orthonormal sketch of $\bQ$}
			\STATE{\;\;\; \makebox[0.5cm]{$\bR$} is $n \times n$ upper triangular R factor}
			\STATE{\textbf{function} $[\bQ, \bS, \bR]=\mathtt{RCholeskyQR}(\bX, \bTheta)$}
			\STATE{1. $\bP \leftarrow \bTheta \bX$}
			\STATE{2. $[\bR,\bS] \leftarrow \mathtt{QR}(\bP)$}
			\STATE{3. $\bQ \leftarrow  \bX \bR^{-1}$}
		\end{algorithmic}
	\end{algorithm}

	\begin{remark}[Sketched SVD]
		In principle, in step 2, we could have orthogonalized $\bP$ with SVD rather than QR factorization. In this case, the matrix $\bR$ would have the form $\bSigma \bV^\mathrm{T}$, where $\bSigma$ is diagonal and $\bV$ is orthonormal. Then step 3 could be stably performed by computing $\bQ \leftarrow (\bX \bV) \bSigma^{-1}$. In this case, multiplying $\bX$ by $\bV$ will require twice as many flops than a forward substitution in  RCholeskyQR. 
		Numerical analysis of such sketched SVD orthogonalization should be aligned with that of RCholeskyQR presented in~\cref{stabRCholeskyQR}.
	\end{remark}

	If the application requires an orthonormal Q factor, and not just a very well-conditioned one, then RCholeskyQR can be augmented with the classical CholeskyQR, resulting in RCholeskyQR2 (see~\cref{alg:RCQR2}). The stability of the RCholeskyQR2 algorithm follows directly from the stability of RCholeskyQR. The computational cost benefit over reorthogonalized shifted CholeskyQR2, i.e., shifted CholeskyQR3, can be characterized in the same way as before. Namely, RCholeskyQR2 should have half the cost in terms of flops, as well as 1.5 times fewer global synchronizations and data passes.
	
	It is important to note that instead of providing the Q factor as an explicit matrix, we can provide it as a function handle that outputs its products with vectors/matrices. Thus, we can overcome the forward substitution in step 4 and therefore reduce complexity of~\cref{alg:RCQR2} by almost a third. Then the left multiplication of the Q factor by given matrix $\bY$ can be performed as $(\bY\bQ){\bR'}^{-1}$, and the right multiplication as $\bQ({\bR' }^{-1} \bY)$. The numerical stability here follows directly from the fact that the matrices $\bQ$ and $\bR'$ are very well conditioned.

	\begin{algorithm}[h] \caption{Augmented randomized Cholesky QR (RCholeskyQR2)} \label{alg:RCQR2}
		\begin{algorithmic}
			\STATE{\textbf{Input:}}
			\STATE{\;\;\; \makebox[0.5cm]{$\bX$} is $m \times n$ matrix} 
			\STATE{\;\;\; \makebox[0.5cm]{$\bTheta$} is $k \times m$ sketching matrix (possibly provided as a function handle)} 		
			\STATE{\textbf{Output}:}
			\STATE{\;\;\; \makebox[0.5cm]{$\bQ$} is $m \times n$ orthonormal Q factor}
			\STATE{\;\;\; \makebox[0.5cm]{$\bR$} is $n \times n$ upper triangular R factor}
			\STATE{\textbf{function} $[\bQ, \bR]=\mathtt{RCholeskyQR2}(\bX,\bTheta)$}
			\STATE{1. $[\bQ, \bS, \bR] \leftarrow \mathtt{RCholeskyQR}(\bX,\bTheta)$} 
			\STATE{2. $\bA \leftarrow \bQ^\mathrm{T} \bQ$}
			\STATE{3. $\bR' \leftarrow \mathtt{chol}(\bA)$, $\bR \leftarrow \bR' \bR$}
			\STATE{4. $\bQ \leftarrow  \bQ \bR'^{-1}$}
		\end{algorithmic}
	\end{algorithm}

	\subsection{Column-oriented RCholeskyQR} \label{colRCholeskyQR}
	Often the columns of $\bX$ are generated recursively from the computed columns of $\bQ$ and $\bR$ in previous iterations.  Let $\bX$ be given by $p$ blocks of columns 
	$$ \bX_{(1:p)} = [\bX_{(1)}, \bX_{(2)}, \hdots, \bX_{(p)} ], $$
	where each block of columns $\bX_{(i)}$ is obtained from the QR factorization $\bQ_{(1:i-1)}\bR_{(1:i-1,1:i-1)}$ of the previously generated matrix $\bX_{(1:i-1)}$. This situation appears, for instance, during the generation of a Krylov basis~\cite{balabanov2021randomized} with Arnoldi iteration:
	$$ \bX_{(i)} \leftarrow \bA \bQ_{(i-1)}, $$
	where $\bA$ is the operator. 
	In such case RCholeskyQR can be performed block by block as shown in~\cref{alg:colRCholeskyQR}.
	\begin{algorithm}[h] \caption{Column-oriented $\mathtt{RCholeskyQR}$} \label{alg:colRCholeskyQR}
		\begin{algorithmic}
			\STATE{\textbf{Input:}}
			\STATE{\;\;\; \makebox[0.6cm]{$\bX_{(1)}$} is $m \times \frac{n}{p}$ matrix} 
			\STATE{\;\;\; \makebox[0.6cm]{$\bTheta$} is $k \times m$ sketching matrix (possibly provided as a function handle)} 		
			\STATE{\textbf{Output}:}
			\STATE{\;\;\; \makebox[0.6cm]{$\bQ$} is $m \times n$ orthonormal Q factor}
			\STATE{\;\;\; \makebox[0.6cm]{$\bS$} is $k \times n$ orthonormal sketch of $\bQ$}
			\STATE{\;\;\; \makebox[0.6cm]{$\bR$} is $n \times n$ upper triangular R factor}
			\STATE{\textbf{function} $[\bQ, \bS, \bR]=\mathtt{colCholeskyQR}(\bX,\bTheta)$}
			\FOR{$i = 1:p$} 
			\STATE{1. If $i>1$ obtain $\bX_{(i)}$ from $\bQ_{(1:i-1)}$ and $\bR_{(1:i-1,1:i-1)}$}
			\STATE{2. $\bP_{(i)} \leftarrow \bTheta \bX_{(i)}$}
			\STATE{3. $\bR_{(1:i-1,i)} \leftarrow \bS_{(1:i-1)}^\dagger \bP_{(i)}$}
			\STATE{4. $[\bS_{(i)}, \bR_{(i,i)}] \leftarrow \mathtt{QR}(\bP_{(i)}-\bS_{(1:i-1)}\bR_{(1:i-1,i)})$}
			\STATE{5.  $\bQ_{(i)} \leftarrow \left  (\bX_{(i)} - \bQ_{(1:i-1)} \bR_{(1:i-1,i)} \right ) \bR_{(i,i)}^{-1}$}
			\ENDFOR
		\end{algorithmic}
	\end{algorithm}
	The least-squares solution in step 3 can be computed with any stable least-squares solver, for instance, based on Householder transformations. The product with $\bR_{(i,i)}^{-1}$ in step 5 of~\cref{alg:colRCQR} is done by forward substitution. It is easy to see that~\cref{alg:colRCholeskyQR} is numerically equivalent to~\cref{alg:RCQR}. Note that when $n = p$,~\cref{alg:colRCholeskyQR} corresponds to the situation when $\bX$ is given column by column.
	
	\Cref{alg:colRCholeskyQR}  has a high relation to the RGS algorithm proposed in~\cite{balabanov2021randomizedGS,balabanov2021randomized}. The main difference is how they update the sketch. In RCholeskyQR algorithm the matrix $\bS_{(i)}$ is effectively computed by a QR factorization of $\bP_{(1:i)}$, while in RGS it is computed as $\bTheta \bQ_{(i)}$. In particular, \Cref{alg:colRCholeskyQR} would exactly recover the RGS algorithm if in step 2 together with computation $\bP_{(i)} \leftarrow \bTheta \bX_{(i)}$ it would also  update the sketch of $\bQ_{(i-1)}$ by calculating $\bS_{(i-1)} \leftarrow \bTheta \bQ_{(i-1)}$.
	Note that this step should increase the cost of the column-oriented RCholeskyQR only by a negligible amount in terms of flops as well as memory consumption and communication cost.
	In general, the RGS algorithm should provide greater stability than RCholeskyQR, since it has the ability to take into account in $\bR_{(:,i)}$ the rounding errors committed when calculating $\bQ$ at previous iterations, while in RCholeskyQR the blocks $\bR_{(:,i)}$ are calculated independently of $\bQ$. This fact is confirmed in our experiments in~\Cref{numexperiments}. Nevertheless, RCholeskyQR shows similar stability as RGS when the input matrix $\bX$ is numerically full-rank, and is slightly less computationally expensive than RGS.

	\subsection{Reduced RCholeskyQR} \label{rRColQR}
	
	Particular attention has to be paid to the case when we are only interested in the low-dimensional projection of the matrix $\bQ$, and not the full matrix. Assume that for given $\bX$ we want to compute the quantity 
	\begin{equation} \label{eq:reduced}
	\bL(\bQ) = \bL \bQ,
	\end{equation}
	where $\bL$ is some (possibly randomized) low-dimensional extractor of the quantity of interest, and $\bQ$ is a well-conditioned matrix satisfying $\mathrm{range}(\bQ) = \mathrm{range}(\bX)$. This situation for instance appears when solving linear system of equation $\bA \bx = \bb$ by a sketched Galerkin or minres projection onto reduced basis $\bX$~\cite{balabanov2019randomized,balabanov2019randomized2}. In such case the $\bL$ extractor has the following form: $$\bL(\bQ) = \bL \bQ = \begin{bmatrix} \bU \bQ \\ \bPhi \bQ  \\ \bPhi (\bA \bQ) \end{bmatrix},$$ where  $\bU$ is (efficient) extractor of low-dimensional quantity $s(\bx) =\bU \bx$ of interest from $\bx$, and  $\bPhi$ is an OSE.  In details, given $\bL(\bQ)$, an approximate solution in the reduced basis can be efficiently and stably obtained by solving the following reduced system of equations
	$$ \bA_{\mathrm{red}} \ba_{\mathrm{red}} = \bb_{\mathrm{red}},$$
	where $\bA_{\mathrm{red}} = (\bPhi \bQ)^\mathrm{T} (\bPhi \bA \bQ)$ and $\bb_{\mathrm{red}} = (\bPhi \bQ)^\mathrm{T} (\bPhi \bb)$ for the sketched Galerkin projection, or $\bA_{\mathrm{red}} = (\bPhi \bA \bQ)^\mathrm{T} (\bPhi \bA \bQ)$ and $\bb_{\mathrm{red}} = (\bPhi \bA \bQ)^\mathrm{T} (\bPhi \bb)$ for the sketched minres projection. Then the linear system's quantity of interest $s(\bx)$ is obtained by calculating $(\bU\bQ) \ba_{\mathrm{red}}$. Furthermore, the above consideration can be extended from the case of solving one linear system to solving a series of systems $\bA(\mu) \bx(\mu) = \bb(\mu)$, with parameters $\mu$ in some set~\cite{balabanov2019randomized,balabanov2019randomized2}.
	
	Then we notice that in order to compute $\bL(\bQ)$, instead of first obtaining $\bQ$ by calculating $\bX\bR^{-1}$ and then applying the extractor $\bL$, we can first apply the extractor to the matrix $\bX$ and only then compute a product with $\bR^{-1}$.
	The resulting reduced RCholeskyQR is depicted in~\cref{alg:oneRCQR}.
	\begin{algorithm}[h] \caption{Reduced $\mathtt{RCholeskyQR}$} \label{alg:oneRCQR}
		\begin{algorithmic}
			\STATE{\textbf{Input:}}
			\STATE{\;\;\; \makebox[0.5cm]{$\bX$} is $m \times n$ matrix} 
			\STATE{\;\;\; \makebox[0.5cm]{$\bTheta$} is $k \times m$ sketching matrix (possibly provided as a function handle)} 
			\STATE{\;\;\; \makebox[0.5cm]{$\bL$} is $l \times m$ extractor of the quantity of interest (possibly provided as a function handle)} 		
			\STATE{\textbf{Output}:}
			\STATE{\;\;\; \makebox[0.5cm]{$\bZ$} is $l \times n$ quantity of interest $\bL(\bQ)$, where $\bQ$ is well-conditioned Q factor}
			\STATE{\;\;\; \makebox[0.5cm]{$\bS$} is $k \times n$ orthonormal sketch of $\bQ$}
			\STATE{\;\;\; \makebox[0.5cm]{$\bR$} is $n \times n$ upper triangular R factor}
			\STATE{\textbf{function} $[\bZ, \bS, \bR ]=\mathtt{redRCholeskyQR}(\bX, \bTheta)$}
			\STATE{1. $\bP \leftarrow \bTheta \bX$, $\bY \leftarrow \bL \bX$}
			\STATE{2. $[\bS, \bR] \leftarrow \mathtt{QR}(\bP)$} 
			\STATE{3.  $\bZ \leftarrow \bY \bR^{-1}$}
		\end{algorithmic}
	\end{algorithm}
	We see that the dominant operation in step~3 of the RCholeskyQR has been drastically reduced. Now the dominant cost comes from the computation of $\bTheta \bX$ and $\bL \bX$ in step 1. By using $\bTheta$ that is SRHT, this operation should have only $\mathcal{O}(m n \log m)$ complexity, which is by a factor $\mathcal{O}(\frac{n}{\log{m}})$ lower than the dominant operations in RCholeskyQR or standard Cholesky QR algorithms. Moreover, it requires only one pass over $\bX$ that can be crucial for the out-of-core computations. \Cref{alg:oneRCQR} has a high relation to the reduced basis orthogonalization depicted in~\cite[Section 4.4]{balabanov2019randomized}.

	Furthermore, if the columns of $\bX$ are generated iteratively from the computed columns of $\bL \bQ$ at the previous iterations, for instance from the reduced linear system's solution $\ba_{\mathrm{red}}$, then the reduced RCholeskyQR can be performed block by block similarly to~\cref{alg:colRCholeskyQR} (see~\cref{alg:colRCQR}).
	\begin{algorithm}[h] \caption{Reduced column-oriented $\mathtt{RCholeskyQR}$} \label{alg:colRCQR}
		\begin{algorithmic}
			\STATE{\textbf{Input:}}
			\STATE{\;\;\; \makebox[0.6cm]{$\bX_{(1)}$} is $m \times \frac{n}{p}$ matrix} 
			\STATE{\;\;\; \makebox[0.6cm]{$\bTheta$} is $k \times m$ sketching matrix (possibly provided as a function handle)} 	
			\STATE{\;\;\; \makebox[0.6cm]{$\bL$} is $l \times m$ extractor of the quantity of interest (possibly provided as a function handle)} 		
			\STATE{\textbf{Output}:}
			\STATE{\;\;\; \makebox[0.6cm]{$\bZ$} is $l \times n$ quantity of interest $\bL(\bQ)$, where $\bQ$ is well-conditioned Q factor}
			\STATE{\;\;\; \makebox[0.6cm]{$\bS$} is $k \times n$ orthonormal sketch of $\bQ$}
			\STATE{\;\;\; \makebox[0.6cm]{$\bR$} is $n \times n$ upper triangular R factor}
			\STATE{\textbf{function} $[\bZ, \bS, \bR ]=\mathtt{colredRCholeskyQR}(\bX, \bTheta)$}
			\FOR{$i = 1:p$} 
			\STATE{1. If $i>1$ obtain $\bX_{(i)}$ from $\bZ_{(1:i-1)}$ and $\bR_{(1:i-1,1:i-1)}$}
			\STATE{2. $\bP_{(i)} \leftarrow \bTheta \bX_{(i)}$, $\bY_{(i)} \leftarrow \bL \bX_{(i)}$}
			\STATE{3. $\bR_{(1:i-1,i)} \leftarrow \bS_{(1:i-1)}^\dagger \bP_{(i)}$}
			\STATE{4. $[\bS_{(i)}, \bR_{(i,i)}] \leftarrow \mathtt{QR}(\bP_{(i)}-\bS_{(1:i-1)}\bR_{(1:i-1,i)})$}
			\STATE{5.  $\bZ_{(i)} \leftarrow \left  (\bY_{(i)} - \bZ_{(1:i-1)} \bR_{(1:i-1,i)} \right ) \bR_{(i,i)}^{-1}$}
			\ENDFOR
		\end{algorithmic}
	\end{algorithm}
	\Cref{alg:colRCQR} is a single-pass algorithm, i.e. it does not require storage/operations with $\bX_{(1:i-1)}$ or $\bQ_{(1:i-1)}$ to get the solution at iteration $i$. Again, it is easy to see that~\cref{alg:oneRCQR} and~\cref{alg:colRCQR} are numerically equivalent. The stability of these algorithms follows from the stability of RCholeskyQR. See~\cref{stabredRCholeskyQR} for more details.

	\section{Rank-revealing randomized Cholesky QR} \label{rrrandCholeskyQR}
	When $\bX$ is numerically rank-deficient, one can improve the stability and the computational cost of RCholeskyQR by orthogonalizing only the linearly independent columns of $\bX$ and ignoring the other columns. In details, we here look for a QR factorization of the form 
	$$ \bX \approx \bQ \bR \bPi^\mathrm{T},$$
	where $\bPi$ is a permutation matrix, $\bQ$ is a very well-conditioned Q factor with $r \leq m$ columns, and $\bR$ is an upper triangular or trapezoidal R factor with $r$ rows. Furthermore, the factorization is such that $\bX \bPi_{(1:r)}$ is sufficiently well-conditioned, and 
	\begin{equation}
	\bX + \Delta \bX = \bQ \bR \bPi^\mathrm{T} \text{ with } \|\Delta \bX_{(:,i)} \|_2 \leq  F(n)u \|\bX_{(:,i)} \|_2
	\end{equation}
	with $F(n)$ being a low-degree polynomial, and $u$ denoting the unit roundoff, $1 \leq i \leq n$. 
	
	To obtain such factorization we propose to first compute $\bPi$ and $\bR$ with the rank-revealing QR of the column-normalized sketch $\bTheta \bX$, and then compute $\bQ$ with forward substitution as is described in~\cref{alg:RRRCholeskyQR}.
	
	\begin{algorithm}[h] \caption{Rank-revealing randomized Cholesky QR ($\mathtt{RRRCholQR}$)} \label{alg:RRRCholeskyQR}
		\begin{algorithmic}
			\STATE{\textbf{Input:}}
			\STATE{\;\;\; \makebox[0.5cm]{$\bX$} is $m \times n$ matrix with normalized columns} 
			\STATE{\;\;\; \makebox[0.5cm]{$\bTheta$} is $k \times m$ sketching matrix (possibly provided as a function handle)} 		
			\STATE{\textbf{Output}:}
			\STATE{\;\;\; \makebox[0.5cm]{$\bQ$} is $m \times r$ well-conditioned Q factor}
			\STATE{\;\;\; \makebox[0.5cm]{$\bS$} is $k \times r$ orthonormal sketch of $\bQ$}
			\STATE{\;\;\; \makebox[0.5cm]{$\bR$} is $r \times n$ upper triangular or trapezoidal R factor}
			\STATE{\;\;\; \makebox[0.5cm]{$\bPi$} is $n \times n$ permutation matrix}
			\STATE{\textbf{function} $[\bQ, \bS, \bR, \bPi]=\mathtt{RRRCholQR}(\bX,\bTheta)$}
			\STATE{1. $\bP \leftarrow \bTheta \bX$}
			\STATE{2. $[\bS, \bR, \bPi] \leftarrow \mathtt{RRQR}(\bP)$} 
			\STATE{3. Determine min $r$ such that $ \|\bR_{(r+1:n,r+1:n)}\|_\Frob \leq {\tau} \|\bR\|_2$}
			\STATE{4.  $\bQ \leftarrow (\bX \bPi_{(:,1:r)}) \bR_{(1:r,1:r)}^{-1}$}
			\STATE{5.  $\bR \leftarrow \bR_{(1:r,:)}$}
		\end{algorithmic}
	\end{algorithm}
	
	For better presentation in~\Cref{alg:RRRCholeskyQR} we assumed that the columns of $\bX$ have unit norms. The case when the columns vary in norm can be accounted for by calculating the normalization matrix $\bD =\mathrm{diag}(\|\bX_{(:, j)}\|_2)$ and inputting $\bX \bD^ { -1}$ to \cref{alg:RRRCholeskyQR} instead of $\bX$. The output R factor has then to be post-processed accordingly to reverse this normalization: $\bR \leftarrow \bR \bPi^\mathrm{T} \bD\bPi$. Furthermore, the computational cost of obtaining $\bD$ in the limited memory architecture and the distributed architecture can be reduced by computing this matrix along with  $\bP \leftarrow \bTheta \bX$ in step 1 during the same pass through the matrix $\bX$ and global synchronization between processors. Then, since we need a sketch of the normalized $\bX$, between steps 1 and 2 we need to normalize $\bP$ by computing $\bP \leftarrow \bP \bD^{-1}$. Also, one can defer the multiplication of $\bX$  by $\bD^{ -1}$ to step 4, which will require operation with only $r$ columns of $\bX$ and hence have a lower computational cost.

	The subroutine $\mathtt{RRQR}$ in step 2 that outputs a rank-revealing QR factorization $\bS\bR\bPi^\mathrm{T}$ of $\bP$ can be chosen for instance  as the strong rank-revealing Cholesky QR from~\cite{gu2004strong} executed in sufficient precision. In this case~\cref{alg:RRRCholeskyQR} is equivalent to the  generalized strong rank-revealing Cholesky QR with respect to the sketched (randomized) inner product, which gives~\cref{alg:RRRCholeskyQR} its name. However, similarly to RCholeskyQR, it should be better to take more robust subroutine in step 2, such as for instance the strong rank-revealing QR from~\cite{gu1996efficient}. 
	
	\Cref{alg:RRRCholeskyQR} takes the truncation tolerance $\tau$ as a user-specified parameter. This parameter can be chosen as $\tau = F(n)u$ to have an approximation of $\bX$ close to machine precision. RRRCholeskyQR can be viewed as an improved version of RCholeskyQR with better stability characteristics. This follows from the fact that $\mathtt{RRQR}$ provides $\bPi$ such that the matrix $\bP \bPi_{(:,1:r)}$ and hence $\bX \bPi_{(:,1:r)}$ is numerically full-rank, given that $\bTheta$ is an OSE of sufficiently large size, and $F(n)$ is sufficiently large. Furthermore, it can be shown that (column-normalized) RRRCholeskyQR factorization is a quasi-optimal rank-$r$ approximation of (column-normalized) $\bX$. These two properties are formalized in~\cref{thm:rrCholeskyQR}. We here restrict ourselves only to the case where $\bTheta$ is an OSE. The analysis for the leverage score sampling matrices is similar.

	\begin{proposition} \label{thm:rrCholeskyQR}
		Consider~\cref{alg:RRRCholeskyQR} with $\mathtt{RRQR}$ in step 2 such that 
		$$ \| \bP  - \bS_{(:,1:r)} \bR_{(1:r,:)} \bPi^\mathrm{T}\|_\Frob \leq C \min_{\mathrm{rank}{(\bY)}=r} \| \bP - \bY \|_\Frob,  $$
		where $C$ is some parameter possibly depending on $n$ and $m$.  
		If $\bTheta$ is an $(\varepsilon, \delta, n)$-OSE, then with probability at least $1-\delta$, we have
		\begin{subequations}
			\begin{equation} \label{eq:lrresult1}
			\sqrt{1-\varepsilon}\| \bX - \bQ \bR  \bPi^\mathrm{T} \|_\Frob  \leq \sqrt{1+\varepsilon} C \min_{\mathrm{rank}{(\bZ)}=r} \| \bX - \bZ \|_\Frob,
			\end{equation}
			and 
			\begin{equation} \label{eq:lrresult2}
			\cond{(\bX  \bPi_{(1:r)})} \leq \sqrt{\frac{1+\varepsilon}{1-\varepsilon}}\cond(\bP \bPi_{(1:r)}). 
			\end{equation}
		\end{subequations}
		Furthermore,~\cref{eq:lrresult1,eq:lrresult2} also hold  with probability at least $1-\delta$, if $\bTheta$ is an $(\varepsilon, (\binom{n}{r+1}+n)^{-1} \delta, r+1)$-OSE and not necessarily an $(\varepsilon, \delta, n)$-OSE. 
		\begin{proof}
			Assume that $r <n$, and that $\bTheta$ is an $\varepsilon$-embedding for all subspaces spanned by $r+1$ columns of $\bX$, and all subspaces of the form $\mathrm{range}(\bZ^{*})+\mathrm{span}(\bX_{(:,j)})$, where $\bX_{(:,j)}$ is the $j$-th column of $\bX$ and $$\bZ^{*} := \arg \min_{\substack{\mathrm{rank}{(\bZ)}=r \\ \mathrm{range}(\bZ) \subset \mathrm{range}(\bX)} } \| \bX - \bZ \|_\Frob,$$  $1\leq j \leq n$. Clearly, this condition is satisfied if $\bTheta$ is an $\varepsilon$-embedding for $\bX$, which in turn holds with probability at least $1-\delta$ if $\bTheta$ is an $(\varepsilon, \delta, n)$-OSE. Furthermore, since there are in total at most $N = \binom{n}{r+1}+n$ such $r+1$-dimensional subspaces, by the union bound argument, $\bTheta$ is an $\varepsilon$-embedding for all of them with probability at least $1-\delta$, if it is an $(\varepsilon, N^{-1}\delta, r+1)$-OSE. 
			
			Notice that since $\bTheta$ is an $\varepsilon$-embedding for all subspaces spanned by $r+1$ columns of $\bX$, it is an $\varepsilon$-embedding for all subspaces of the form $\mathrm{range}(\bQ)+\mathrm{span}(\bX_{(:,j)})$.	
			It then follows that 
			\begin{align*}
			&\| \bP \bPi - \bS_{(:,1:r)} \bR_{(1:r,:)} \|^2_\Frob = \| \bTheta(\bX \bPi - \bX \bPi_{(1:r)} {\bR_{(1:r,1:r)}}^{-1} \bR_{(1:r,:)}) \|^2_\Frob = \sum^n_{i=1} \| \bTheta (\bX_{(:,i)}  - \bQ \bR_{(1:r,:)} \bPi_{(i,:)}^\mathrm{T}) \|^2_2 \\
			& \geq (1-\varepsilon) \sum^n_{i=1} \| \bX_{(:,i)} - \bQ \bR_{(1:r,:)}\bPi_{(i,:)}^\mathrm{T} \|^2_2 =  (1-\varepsilon) \| \bX  - \bQ \bR \bPi^\mathrm{T} \|^2_\Frob,
			\end{align*}
			and 
			\begin{align*}
			\min_{\mathrm{rank}{(\bY)}=r} \| \bP - \bY\|^2_\Frob \leq  \| \bTheta(\bX - \bZ^*)\|^2_\Frob \leq (1+\varepsilon) \sum^n_{i=1} \| \bX_{(:,i)} - \bZ^*_{(:,i)} \|^2_2 = (1+\varepsilon) \| \bX - \bZ^* \|^2_\Frob , 
			\end{align*}
			which gives~\cref{eq:lrresult1}.
			
			The relation~\cref{eq:lrresult2} follows from the fact that for any vector $\bx \in \mathbb{R}^{r}$ we have
			$$  (1-\varepsilon) \| \bX \bPi_{(1:r)} \bx \|^2_2 \leq \|\bTheta \bX \bPi_{(1:r)} \bx \|^2_2 \leq (1+\varepsilon) \| \bX \bPi_{(1:r)} \bx \|^2_2. $$
		\end{proof}
	\end{proposition}
	
	\Cref{thm:rrCholeskyQR} opens the door to yet another application of the RRRCholeskyQR algorithm beyond QR factorization of tall-and-skinny matrices, which is an efficient and stable low-rank approximation of large $m \times n$ matrices that do not need to satisfy $n \ll m$.
	According to~\cref{thm:rrCholeskyQR},~\cref{alg:RRRCholeskyQR} provides a rank-revealing QR factorization of the same quality (in exact arithmetic) as the low-dimensional $\mathtt{RRQR}$ in step 2, given that $\bTheta$ is a $(\varepsilon, \delta, n)$-OSE or $(\varepsilon, (\binom{n}{r+1}+n)^{-1} \delta, r +1 )$-OSE. The first-mentioned condition on $\bTheta$ is met if $\bTheta$ is a Gaussian matrix with $k = \mathcal{O}(n+\log{\frac{1}{\delta}})$ rows, which is sufficient when $\bX$ is tall-and-skinny, but infeasible when both dimensions of $\bX$ are large. In the latter case, we must turn to the second condition, which is satisfied by a Gaussian OSE with $k=\mathcal{O}(r\log{n}+\log\frac{1}{\delta})$ rows. For SRHT this requirement is higher but in practice there should not be much difference in accuracy.  It is concluded that RRRCholeskyQR can provide a quasi-optimal low-rank approximation using the sketching dimension $k$ that depends on $n$ and $m$ at most logarithmically. 
	
	The numerical stability of RRRCholeskyQR in both the low-rank approximation context as well as the context of factorization of tall-and-skinny matrices, can be characterized in exactly the same manner. It can be guaranteed unconditionally of $\bX$. In particular, it can be shown that by using $\tau \geq G(n) u$, where $G(n)$ is some low-degree polynomial, we obtain  $\bP \bPi_{(1:r)}$ and hence $\bX \bPi_{(1:r)}$ of numerically full rank. In turn this fact implies the stability of computing the factor $\bQ$, and therefore the overall stability of the algorithm. More details on this matter are provided in~\cref{stabRRRCholeskyQR}. Furthermore, if the application requires an orthonormal to machnie precision Q factor  and not just well-conditioned, then RRRCholeskyQR can be augmented with the classical CholeskyQR in exactly the same way as done in RCholeskyQR2 defined by~\cref{alg:RCQR2}. The resulting algorithm will be referred to as RRRCholeskyQR2.
	
	In terms of efficiency, RRRCholeskyQR and RRRCholeskyQR2 should be at least as good as RCholeskyQR and RCholeskyQR2 respectively, and therefore outperform standard/shifted CholeskyQR2 and CholeskyQR3, Householder QR and other deterministic algorithms. In addition, it is revealed that when $\bX$ has a relatively low numerical rank, the proposed factorizations can be even more computationally advantageous. In particular, if we take $\bTheta$ as an SRHT matrix, then RRRCholeskyQR should take about $mn \log_2 m+mr^2$ flops in steps 1 and 4, whereas RRRCholeskyQR2 (in implicit form) takes about $mn \log_2 m+2mr^2$ flops.  In contrast, the standard QR factorizations and RCholeskyQR consume $\mathcal{O}(mn^2)$ flops, which can be much larger when the rank $r$ is relatively small. Furthermore, RRRCholeskyQR and RRRCholeskyQR2 should significantly outperform other randomized rank-revealing QR factorizations such as the ones from~\cite{martinsson2015blocked,martinsson2017householder,xiao2017fast,duersch2017randomized,mary2015performance} that require at least $mnr$ flops (needed to compute the R factor) and two or more global synchronizations between processors. The RRRCholeskyQR factorization on the other hand needs up to $\frac{n}{r}$ less flops and only one global synchronization.

	\section{Stability analysis} \label{stability}
	In this section we analyze numerical stability of~\cref{alg:RCQR,alg:oneRCQR,alg:RRRCholeskyQR}. The stability of other presented algorithms directly follows. Let $\kappa$ denote the condition number of column-normalized $\bX$.
	Let us recall that the steps requiring a minor computational cost are here executed in higher precision by a low-degree polynomial factor $F(m,n)$ in $n$ and $m$ than the dominant operations. This allows to have numerical stability with working unit roundoff $u = \mathcal{O}(n^{-\frac{3}{2}} \kappa^{-1})$ in dominant step 3 of~\cref{alg:RCQR,alg:oneRCQR} and $u = \mathcal{O}(n^{-\frac{3}{2}} r^{-\frac{5}{2}})$ in dominant step 4 of~\cref{alg:RRRCholeskyQR} independent of the dimension $m$. 
	Nevertheless, clearly our results also imply guarantees of numerical stability under a model with unique unit roundoff. Such guarantees can be obtained simply by replacing $u$ by $F(m,n)u$ in the results.
	
	Note that the above conditions on $u$ are pessimistic. This overestimation can be seen as an artifact due to the use of worst-case rounding bounds. According to the ``rule of thumb'' of rounding~\cite{higham2002accuracy}, in practice the low-dimensional polynomials in the forthcoming in this section    conditions~\cref{eq:ass0,eq:ass30,eq:ass20} on $u$, and the stability guarantees in~\cref{thm:main,thm:rmain,thm:RRRCholeskyQR1,thm:RRRCholeskyQR2} can be reduced by (nearly) a square root.
	
	\subsection{Stability of RCholeskyQR} \label{stabRCholeskyQR}
	The stability analysis of~\cref{alg:RCQR} is carried out using the following assumptions. First, it is assumed that the working unit roundoff $u$ satisfies the bound~\cref{eq:ass0}. 
	Furthermore, steps 1 and 2 are assumed to be performed with unit roundoff $u_f$ such that the error matrices 
	\begin{subequations} \label{eq:errmat1}
		\begin{align}
		\bE_1  &:= \bTheta \bX - \bhP \\
		\bE_2  &:= \bhP - \bhS \bhR \\
		\bE_3  &:= \bX - \bhQ \bhR.
		\end{align}
	\end{subequations}
	satisfy~\cref{eq:ass1,eq:ass2,eq:ass3}. By the classical worst-case rounding analysis we have 
	$|\bE_1| \leq \frac{m u_f}{1-m u_f}  | \bTheta ||\bX|. $
	Hence~\cref{eq:ass1} can be achieved with $u_{f} = \mathcal{O}(m^{-1} n^{-1} u)$.
	In principle the guarantees~\cref{eq:ass2} can be achieved with any stable QR factorization executed in sufficient precision. This includes the Householder QR with unit roundoff $u_{f}=O(k^{-1}n^{-\frac{3}{2}}u)$ or Givens QR with unit roundoff $u_{f}=O(k^{-1}n^{-\frac{1}{2}}u)$~\cite{higham2002accuracy}.
	According to \cite[Theorem 8.5]{higham2002accuracy}, we have $\bhQ_{(j,:)} (\bhR + \Delta \bR^{(j)}) = \bX_{(j,:)} $ with $|\Delta \bR^{(i)}| \leq 1.1 u n |\bhR|$, which implies that $|\bE_3|$ is bounded by $1.1 u n |\bhQ||\bhR |$ and leads to the first inequality in~\cref{eq:ass3}. Since the rows of $\bE_3$ are computed independently of each other, we can use~\cref{thm:thetadeltax2} to bound the sketched norms of the columns of $\bE_3$. In this way we have $\| \bTheta {\bE_3}_{(:,j)}\|_2 \leq \sqrt{\frac{3}{2}} 1.1 u n \||\bhQ| |\bhR_{(:,j)} | \|_2$, with probability at least $1-\delta$, if $\bTheta$ is $(1/8,n^{-1}\binom{m}{d}^{-1}\delta, d)$-OSE, with $d = 4.2 c^{-1} \log{{\textstyle \frac{4}{\delta}}}$, which in turn is satisfied by Gaussian matrices and SRHT (in practice) with $k \geq \mathcal{O}(\log m \log{\frac{1}{\delta}})$ rows. Note that the bound for 
	$\| \bTheta {\bE_3}_{(:,j)}\|_2$ is independent of the high dimension $m$.

	\begin{assumptions} \label{ass:assumptions1}
		Consider~\cref{alg:RCQR}. 	
		We assume that
		\begin{subequations}
			\begin{equation} \label{eq:ass0} 
			u \leq 0.01 n^{-\frac{3}{2}} \kappa^{-1}.
			\end{equation}
			Furthermore, for $1\leq j \leq n$,
			\begin{align}
			\|{\bE_1}_{(:,j)}\|_2 &\leq 0.1u n^{-\frac{1}{2}} \| \bX_{(:,j)}\|_2, ~~~~~~~~~~~~~~~~~~~~~~~~ \label{eq:ass1} \\
			\|{\bE_2}_{(:,j)}\|_2 &\leq 0.1 u n^{-\frac{1}{2}}  \| \bhP_{(:,j)} \|_2 , ~~~~~~~~~~~~ &\|\bhS^\mathrm{T} \bhS - \bI \|_\Frob \leq 0.1 u,~~~~~~~~~~~~~~~~~~~~~~~~~~~  \label{eq:ass2}\\
			|{\bE_3}_{(:,j)}| &\leq 1.1 u n  | \bhQ| | \bhR_{(:,j)}| ,~~~ &\|\bTheta {\bE_3}_{(:,j)}\|_2 \leq 2 u n \| \bhQ\|_\Frob \| \bhR_{(:,j)}\|_2.  ~~~~~~~~~ \label{eq:ass3}	
			\end{align}
		\end{subequations}		
	\end{assumptions}
	
	\Cref{thm:main} provides a stability guarantee of~\cref{alg:RCQR}.
	
	\begin{theorem} \label{thm:main}
		Let $\bTheta$ be an $\varepsilon$-embedding for $\bX$ with $\varepsilon \leq \frac{1}{2}$.
		Consider~\cref{alg:RCQR}. Under~\cref{ass:assumptions1} we have,
		\begin{subequations}
			\begin{align}	
			\bX +\Delta \bX = \bhQ \bhR \text{ with }& \|\Delta \bX_{(:,j)}\|_2\leq  2.1 u n \|\bX_{(:,j)} \|_2  \label{eq:main1} \\
			(1+\varepsilon)^{-\frac{1}{2}} - 4 u n^\frac{3}{2}  \kappa \leq \sigma_{min}(\bhQ)&\leq \sigma_{max}(\bhQ)  \leq (1-\varepsilon)^{-\frac{1}{2}}+ 4 u n^\frac{3}{2} \kappa,\label{eq:main2}
			\end{align}
			for $1\leq j \leq n$. Furthermore, it holds that
			\begin{equation} \label{eq:main3}
			\|\bhS - \bTheta \bhQ \|_\Frob \leq 6.1 u n^\frac{3}{2}  \kappa,
			\end{equation}
		\end{subequations}
		and that
		$\bTheta$ satisfies the $\varepsilon'$-embedding property for $\bhQ$ with $\varepsilon' \leq \varepsilon +  50 u n^\frac{3}{2}  \kappa$.
		\begin{proof}
			Let us scale $\bX$, $\bhR$, $\bhP$  and $\bE_1$, $\bE_2$, $\bE_3$ by $\bD:=\mathrm{diag}(\|\bX_{(:,j)}\|^{-1}_2)$:  $\bX \leftarrow \bX \bD$, $\bhR \leftarrow \bhR \bD$, $\bhP \leftarrow \bhP \bD$, $\bE_1\leftarrow \bE_1 \bD$, $\bE_2\leftarrow \bE_2 \bD$, and $\bE_3\leftarrow \bE_3 \bD$. Notice that such scaling does not affect the relations~\cref{eq:errmat1} and assumptions~\cref{eq:ass1,eq:ass2,eq:ass3}. Then we also have $\kappa = \mathrm{cond}(\bX)$.

			We start with showing that the computed sketch $\bhP$ of $\bX$ preserves the column norms and the smallest singular value of $\bX$. By the $\varepsilon$-embedding property of $\bTheta$ and~\cref{eq:ass1}, we get for $1 \leq j \leq n$,
			\begin{subequations} \label{eq:stab1}
				\begin{align}
				\|\bhP_{(:,j)} \|_2 &\leq \|\bTheta \bX_{(:,j)} \|_2+\|{\bE_1}_{(:,j)} \|_2 \leq (\sqrt{1+\varepsilon}+0.1un^{-\frac{1}{2}}) \|\bX_{(:,j)} \|_2 \leq 1.23 \|\bX_{(:,j)} \|_2, \\
				\sigma_{min}(\bhP) &\geq \sigma_{min}(\bTheta \bX)-\|\bE_1 \|_2 \geq \sqrt{1-\varepsilon}\sigma_{min}(\bX) - 0.1 u n^{-\frac{1}{2}} \mathrm{cond}(\bX) \sigma_{min}(\bX) \geq 0.69 \sigma_{min}(\bX).
				\end{align}
			\end{subequations}
			Next, the same thing is shown for the computed R factor. We have,
			\begin{equation*}
			\bhR_{(:,j)} = \bhS^\mathrm{T} \bhP_{(:,j)}	+ (\bI-\bhS^\mathrm{T} \bhS)\bhR_{(:,j)} + \bhS^\mathrm{T} {\bE_2}_{(:,j)}.
			\end{equation*}
			Hence by~\cref{eq:ass2,eq:stab1}, it holds
			\begin{subequations} \label{eq:stab2}
				\begin{align}
				\|\bhR_{(:,j)}\|_2  &\leq \|\bhS\|_2 \|\bhP_{(:,j)}\|_2 + \|\bI-\bhS^\mathrm{T} \bhS\|_2\|\bhR_{(:,j)}\|_2 + \|\bhS\|_2 \|{\bE_2}_{(:,j)}\|_2 \leq 1.26 \|\bX_{(:,j)}\|_2,	 \\			
				\sigma_{min}(\bhR)  &\geq \sigma_{min}(\bhS) \sigma_{min} (\bhP) - \|\bI-\bhS^\mathrm{T} \bhS\|_\Frob\|\bhR\|_2 - \|\bhS\|_2 \|\bE_2\|_\Frob \geq  0.66 \sigma_{min}(\bX). 
				\end{align}
			\end{subequations}
			
			This allows us to bound the error of $\bhQ$ in the Frobenius norm and the sketched Frobenius norm. By~\cref{eq:ass3,eq:stab2}, we get
			\begin{subequations} \label{eq:stab3}
				\begin{align}
				\|\bhQ - \bX \bhR^{-1}\|_\Frob  &= \|\bE_3 \bhR^{-1}\|_\Frob \leq \|\bE_3 \|_\Frob \|\bhR^{-1}\|_2 \leq 1.1 u n \| \bhQ\|_2 \| \bhR\|_\Frob \|\bhR^{-1}\|_2 \leq 2.1 u n^\frac{3}{2} \mathrm{cond}(\bX) \| \bhQ\|_2, \\ 
				\|\bTheta\bhQ - \bTheta\bX \bhR^{-1}\|_\Frob  &= \|\bTheta \bE_3 \bhR^{-1}\|_\Frob \leq \|\bTheta\bE_3 \|_\Frob \|\bhR^{-1}\| \leq 2 u n \| \bhQ\|_2 \| \bhR\|_\Frob \|\bhR^{-1}\|_2 \leq 3.82  u n^\frac{3}{2} \mathrm{cond}(\bX) \| \bhQ\|_2
				\end{align}
			\end{subequations}
			
			Furthermore, from~\cref{eq:ass1,eq:ass2,eq:stab2} we obtain
			\begin{equation} \label{eq:stab4}
			\|\bTheta\bX \bhR^{-1} -\bhS\|_\Frob  = \| \bE_1 \bhR^{-1} + \bE_2\bhR^{-1}\|_\Frob \leq 0.1u  n^{-\frac{1}{2}} \|\bX \|_\Frob \|\bhR^{-1}\|_2 + 0.1u n^{-\frac{1}{2}}  \|\bhP\|_\Frob \|\bhR^{-1}\|_2 \leq 0.35  u \mathrm{cond}(\bX),
			\end{equation}
			which in turn implies that 
			\begin{equation} \label{eq:stab5}
			1-\Delta_1 \leq \sigma_{min}(\bTheta \bX \bhR^{-1}) \leq \sigma_{max}(\bTheta \bX \bhR^{-1}) \leq 1+\Delta_1,
			\end{equation}
			with $\Delta_1 \leq \|\bhS^\mathrm{T} \bhS - \bI \|_\Frob + 0.35  u \mathrm{cond}(\bX) \leq  0.45u \mathrm{cond}(\bX)$. By the $\varepsilon$-embedding property of $\bTheta$, it is deduced from~\cref{eq:stab5} that
			\begin{equation}\label{eq:stab6}
			(1+\varepsilon)^{-\frac{1}{2}}(1-\Delta_1) \leq \sigma_{min}(\bX \bhR^{-1}) \leq \sigma_{max}(\bX \bhR^{-1}) \leq (1-\varepsilon)^{-\frac{1}{2}}(1+\Delta_1),
			\end{equation}
			which, combined with~\cref{eq:stab3}, results in
			\begin{equation}\label{eq:stab7}
			(1+\varepsilon)^{-\frac{1}{2}}(1-\Delta_1)-\Delta_2 \leq \sigma_{min}(\bhQ) \leq \sigma_{max}(\bhQ) \leq (1-\varepsilon)^{-\frac{1}{2}}(1+\Delta_1)+\Delta_2,
			\end{equation}
			for some $\Delta_2 \leq 2.1un^\frac{3}{2} \mathrm{cond}(\bX) \|\bhQ\|_2$.
			Hence $\|\bhQ\|_2 \leq 1.5$, which combined with~\cref{eq:ass3,eq:stab2,eq:stab7}, implies the results~\cref{eq:main1,eq:main2} of the theorem. 
			
			Relations~\cref{eq:stab3,eq:stab4} imply~\cref{eq:main3}. To show that $\bTheta$ is an $\varepsilon'$-embedding for $\bhQ$, we notice that for any $\bx \in \mathbb{R}^k$, it holds 
			\begin{equation}\label{eq:stab8}
			(1 - \Delta_3) \| \bx \|_2 \leq \|\bTheta \bX \bhR^{-1} \bx \|_2 - \|\bTheta \bE_3 \bhR^{-1}\bx\|_2 \leq \|\bTheta \bhQ \bx\|_2 \leq \|\bTheta \bX \bhR^{-1} \bx \|_2 + \|\bTheta \bE_3 \bhR^{-1}\bx\|_2 \leq (1 + \Delta_3) \| \bx \|_2,
			\end{equation}
			and
			\begin{equation}\label{eq:stab9}
			\left ((1+ \varepsilon)^{- \frac{1}{2}} - \Delta_3 \right ) \| \bx \|_2\leq \| \bhQ \bx \|_2 \leq  \left ((1- \varepsilon)^{- \frac{1}{2}} + \Delta_3 \right ) \| \bx \|_2, 
			\end{equation}
			where $\Delta_3 \leq 6.1 u n^\frac{3}{2}  \mathrm{cond}(\bX)$. 	
			Whence, 	
			$$( 1 - \Delta_3  )^2 \left ((1- \varepsilon)^{- \frac{1}{2}} + \Delta_3 \right )^{-2}  \|\bhQ \bx\|^2_2 \leq  \|\bTheta \bhQ \bx\|^2_2 \leq   ( 1 + \Delta_3  )^2 \left ((1+ \varepsilon)^{- \frac{1}{2}} - \Delta_3 \right )^{-2}  \|\bhQ \bx\|^2_2. $$
			By using the fact that $\Delta_3 \leq 0.061$, we get
			$( 1 + \Delta_3  )^2 \left ((1+ \varepsilon)^{-\frac{1}{2}} - \Delta_3 \right )^{-2} = (1+\varepsilon) ( 1 + \Delta_3  )^2  (1-\sqrt{1+\varepsilon}\Delta_3)^{-2} \leq 1+\varepsilon+8 \Delta_3,	$  and similarly,
			$( 1 - \Delta_3  )^2 \left ((1- \varepsilon)^{-\frac{1}{2}} + \Delta_3 \right )^{-2}  \geq 1-\varepsilon-8 \Delta_3.	$
			
		\end{proof}
	\end{theorem}

	The stability of the column-oriented RCholeskyQR given by~\cref{alg:colRCholeskyQR} follows directly from~\cref{thm:main}, since it is numerically equivalent to~\cref{alg:RCQR}.

	\subsection{Stability of reduced RCholeskyQR} \label{stabredRCholeskyQR}
	
	Stability guarantees for the reduced RCholeskyQR and its column-oriented variant (\Cref{alg:colRCQR,alg:oneRCQR}) can be obtained in a similar manner as above. For consistency with the previous subsection, we again assume that the forward substitution in step 3 in~\cref{alg:oneRCQR} and in step 5 in~\Cref{alg:colRCQR} is done with unit roundoff $u$, which is by a polynomial factor in $n$ and $m$ less than the unit rounding $u_f$ used for other operations. Although for a reduced RCholeskyQR, this condition may not be that important, since the forward substitution is performed on a low-dimensional matrix and, therefore, should not be so expensive.
	
	We shall analyze only~\Cref{alg:oneRCQR}, noting that the results obtained will also apply to~\Cref{,alg:colRCQR}, since the two algorithms are numerically equivalent.  Let $\bL$ have $l \ll m$ rows. Let $\bE_1$ and $\bE_2$ be the rounding matrices defined in~\cref{eq:errmat1}. Also, define
	\begin{subequations} \label{eq:errmat3}
		\begin{align}
		\bE_4 &:= \bhY - \bL \bX \\
		\bE_5 &:= \bhZ \bhR - \bhY. 
		\end{align}
	\end{subequations}
	Our analysis will be based on the following assumptions. 
	\begin{assumptions} \label{ass:assumptions3}
		Consider~\cref{alg:colRCQR}. We assume that 
		\begin{subequations}
			\begin{equation} \label{eq:ass30} 
			u \leq 0.01 n^{-\frac{3}{2}} \kappa^{-1}.
			\end{equation}
			Furthermore, for $1 \leq i \leq l$, $1\leq j \leq n$,
			\begin{align}
			\|{\bE_1}_{(:,j)}\|_2 &\leq 0.1u n^{-\frac{1}{2}} \| \bX_{(:,j)}\|_2, ~~~~~~~~~~~~~~~~~~~~~~~~ \label{eq:ass31} \\
			\|{\bE_2}_{(:,j)}\|_2 &\leq 0.1 u n^{-\frac{1}{2}}  \| \bhP_{(:,j)} \|_2 , ~~~~~~~~~~~~ &\|\bhS^\mathrm{T} \bhS - \bI \|_\Frob \leq 0.1 u,~~~~~~~~~~~~  \label{eq:ass32}\\
			|{\bE_4}_{(i,j)}| &\leq 0.1 u  \|\bL_{(i,:)}\|_2 \|\bX_{(:,j)}\|_2, \label{eq:ass33} \\
			|{\bE_5}_{(i,j)}| &\leq 1.1 u n  \| \bhZ_{(i,:)}\|_2 \| \bhR_{(:,j)}\|_2. \label{eq:ass34}
			\end{align}
		\end{subequations}
	\end{assumptions}
	
	The conditions~\cref{eq:ass31,eq:ass32} are the same as~\cref{eq:ass1,eq:ass2} in the analysis of RCholeskyQR. They can be achieved by using unit roundoff $u_{f} = \mathcal{O}(m^{-1} n^{-1} u)$. With this unit roundoff we also can get~\cref{eq:ass33} that follows directly from standard rounding analysis from~\cite{higham2002accuracy}. Furthermore, \cite[Theorem 8.5]{higham2002accuracy} states that $\bhZ_{(i,:)} (\bhR + \Delta \bR^{(i)}) = \bhY_{(i,:)} $ with $|\Delta \bR^{(i)}| \leq 1.1 u n |\bhR|$ that in turn implies~\cref{eq:ass34}.

	\Cref{thm:rmain} provides stability characterization of~\cref{alg:oneRCQR}.
	\begin{theorem} \label{thm:rmain}
		Let $\bTheta$ be an $\varepsilon$-embedding for $\bX$ with $\varepsilon \leq \frac{1}{2}$. 	Consider~\cref{alg:oneRCQR}. Under~\cref{ass:assumptions3}, there exists $\bQ$ with $\mathrm{range}(\bQ) = \mathrm{range}(\bX)$ such that
		\begin{subequations}
			\begin{align}	
			\| \bhZ_{(i,:)} - \bL_{(i,:)} \bQ\|_2 &\leq 3.5 u n^{\frac{3}{2}} \kappa \| \bL_{(i,:)} \|_2,  \label{eq:rmain1}  \\
			(1+\varepsilon)^{-\frac{1}{2}} - 0.45 u  \kappa \leq \sigma_{min}(\bQ)&\leq \sigma_{max}(\bQ)  \leq (1-\varepsilon)^{-\frac{1}{2}}+ 0.45 u \kappa, \label{eq:rmain2}
			\end{align}
			for $1 \leq i \leq l$.
			In addition,
			\begin{equation} \label{eq:rmain3}
			\|\bhS - \bTheta \bQ \|_\Frob \leq 0.35 u \kappa
			\end{equation}
		\end{subequations}
		and
		$\bTheta$ is an $\varepsilon$-embedding for $\bQ$.
		
		\begin{proof}
			Take $ \bQ = \bX \bhR^{-1}.$ As in~\cref{thm:main}, scale $\bX$, $\bhR$, $\bhP$, $\bhY$  and $\bE_1$, $\bE_2$, $\bE_4$, $\bE_5$ by $\bD=\mathrm{diag}(\|\bX_{(:,j)}\|^{-1}_2)$:  $\bX \leftarrow \bX \bD$, $\bhR \leftarrow \bhR \bD$, $\bhP \leftarrow \bhP \bD$, $\bhY \leftarrow \bhY \bD$, $\bE_1\leftarrow \bE_1 \bD$, $\bE_2\leftarrow \bE_2 \bD$, $\bE_4\leftarrow \bE_4 \bD$, $\bE_5\leftarrow \bE_5 \bD$, which does not affect the assumptions.
			
			The fact that $\bTheta$ is an $\varepsilon$-embedding for $\bQ$ is obvious. Furthermore, the results~\cref{eq:rmain2,eq:rmain3} can be proven similarly to~\cref{eq:stab4,eq:stab5} in the proof of~\cref{thm:main}.	
			
			To show~\cref{eq:rmain1} we shall use the following result (see~\cref{eq:stab2}) from the proof of~\cref{thm:main}:			
			\begin{align}		 
			\|\bhR \|_\Frob \leq 1.26 \|\bX \|_\Frob~&\text{ and } \sigma_{min}(\bhR) \geq 0.66 \sigma_{min}(\bX). \label{eq:stab2r} 
			\end{align}	
			Then by~\cref{eq:ass33,eq:ass34} we have		
			
			\begin{align*}
			\| \bZ_{(i,:)} - \bL_{(i,:)} \bQ\|_2 &\leq  \| \bZ_{(i,:)} - \bY_{(i,:)} \bhR^{-1}\|_2 + \|\bY_{(i,:)} \bhR^{-1} - \bL_{(i,:)} \bX \bhR^{-1} \|_2 \\
			&\leq \| {\bE_5}_{(i,:)}\|_2 \|\bhR^{-1}\|_2 + \|{\bE_4}_{(i,:)}\|_2\|\bhR^{-1} \|_2 \\
			&\leq 1.1 u n^{\frac{3}{2}} \| \bZ_{(i,:)} \|_2    \| \bhR \|_2 \|\bhR^{-1} \|_2 + 0.1 u n^{\frac{1}{2}}  \| \bL_{(i,:)} \|_2 \| \bX \|_2 \|\bhR^{-1} \|_2 \\
			& \leq 1.91 u n^{\frac{3}{2}} \mathrm{cond}(\bX) \| \bZ_{(i,:)} \|_2 + 0.2 u n^{\frac{1}{2}} \mathrm{cond}(\bX)  \| \bL_{(i,:)} \|_2.
			\end{align*}
			This relation particularly implies that $\|\bZ_{(i,:)} \|_2 \leq 1.6 \|\bL_{(i,:)} \|_2.$ Consequently, we have
			$$ \| \bZ_{(i,:)} - \bL_{(i,:)} \bQ\|_2 \leq 3.5 u n^{\frac{3}{2}} \mathrm{cond}(\bX) \|\bL_{(i,:)} \|_2, $$
			that is equivalent to~\cref{eq:rmain1}.
		\end{proof}
	\end{theorem}
	According to~\cref{thm:rmain}, \Cref{alg:oneRCQR} computes a quantity of interest  $\bZ = \bL \bQ$ associated with some well-conditioned matrix $\bQ$ with $\mathrm{range}(\bQ) = \mathrm{range}(\bX)$, with relative row-wise errors 
	$$ \frac{\| \bhZ_{(i,:)} - \bL_{(i,:)} \bQ\|_2}{\| \bL_{(i,:)} \bQ \|_2} = \mathcal{O}(u  n^{\frac{3}{2}} \kappa),~~~1 \leq i \leq l.$$

	\subsection{Stability of RRRCholeskyQR} \label{stabRRRCholeskyQR}
	This section is devoted to stability analysis of the RRRCholeskyQR algorithm. To simplify the presentation, we redefine $\bX$ by permuting its columns with the permutation matrix $\bPi$:
	$$\bX \leftarrow \bX \bPi.$$
	In addition, it is assumed that $\bX$ has normalized columns, since the errors due to normalization are here negligible and can be ignored.
	
	First notice that the factorization $\bX_{(1:r)} = \bhQ \bhR_{(1:r)}$ can be viewed as a RCholeskyQR factorization of $\bX_{(1:r)}$. Define the associated error matrices: 
	\begin{subequations} \label{eq:errmat2}
		\begin{align}
		\bE^{*}_1  &:= \bTheta \bX_{(1:r)} - \bhP_{(1:r)}, \\
		\bE^{*}_2  &:= \bhP_{(1:r)} - \bhS \bhR_{(1:r)}, \\
		\bE^{*}_3  &:= \bX_{(1:r)} - \bhQ \bhR_{(1:r)}.
		\end{align}
	\end{subequations}
	Define also the error matrix associated with the computation of $\bTheta \bX_{(r+1:n)}$:
	\begin{equation} 
	\bE_6  := \bTheta \bX_{(r+1:n)} - \bhP_{(r+1:n)}.
	\end{equation}
	The stability analysis will be based on the following assumptions. 
	\begin{assumptions} \label{ass:ass200}
		Consider~\cref{alg:RRRCholeskyQR}. We assume that
		\begin{subequations} \label{eq:ass200} 
			\begin{equation} \label{eq:ass20} 
			u \leq 0.001 n^{-\frac{3}{2}} r^{-\frac{5}{2}}.
			\end{equation}
			Furthermore,
			\begin{align}
			&\|{\bE^{*}_1}_{(:,j)}\|_2 \leq 0.01ur^{-\frac{1}{2}}   \| \bX_{(:,j)} \|_2,  & ~~~~~~~~~~~~~~~~~~~~~~~~ \label{eq:ass21} \\
			&\|{\bE^{*}_2}_{(:,j)}\|_2 \leq 0.01  ur^{-\frac{1}{2}}   \| \bP_{(:,j)} \|_2 , ~~~~~~~~~~~~ &\|\bhS^\mathrm{T} \bhS - \bI \|_\Frob \leq 0.1 u,~~~~~~~~~~~~  \label{eq:ass22}\\
			&|{\bE^{*}_3}_{(:,j)}| \leq 1.1 u r  | \bhQ| | \bhR_{{(:,j)}}| ,~~~ &\|{\bTheta \bE^{*}_3}_{(:,j)}\|_\Frob \leq 2 u  r \| \bhQ\|_\Frob \| {\bhR}_{(:,j)}\|_2 ,  ~~~~~~~~~~~~~ \label{eq:ass23}	\\
			&\|{\bE_6}_{(:,j)}\|_2 \leq 0.1u n^{-\frac{1}{2}}\| \bX_{(:,j)} \|_2.  & ~~~~~~~~~~~~~~~~~~~~~~~~ \label{eq:ass24}
			\end{align}
			We also assume that the $\mathtt{RRQR}$ subroutine in step 2 is such that 
			\begin{align}
			&\sigma_r(\bhR_{(1:r)}) \geq 0.5 n^{-\frac{1}{2}} r^{-\frac{1}{2}} \sigma_r(\bhP),  &\|\bhR_{(r:n)}\|_2 \leq 2n^\frac{1}{2} r^\frac{1}{2} \sigma_{r}(\bhP), ~~~~~~~~~~~~~~~~~~~~~~~~ \label{eq:ass25}	\\
			&\|\bhR_{(1:r)}^{-1} \bhR_{(r+1:n)} \|_\Frob  \leq 2 n^\frac{1}{2} r^\frac{1}{2}. &  ~~~~~~~~~~~~~~~~~~~~~~~~\label{eq:ass26}	
			\end{align}	
		\end{subequations}
	\end{assumptions}

	The assumptions~\cref{eq:ass21,eq:ass22,eq:ass23} repeat~\cref{eq:ass1,eq:ass2,eq:ass3} when $\bX_{(1:r)} = \bhQ \bhR_{(1:r)}$ is seen as a RCholeskyQR factorization of $\bX_{(1:r)}$. They can be satisfied by using unit roundoff $u_{f} = \mathcal{O}(m^{-1} n^{-1} u)$ in the minor steps of~\cref{alg:RRRCholeskyQR}. Furthermore, by the standard rounding analysis we have $|\bE_6|\leq \frac{m u_f}{1-m u_f}  | \bTheta ||\bX_{(r+1:n)}|, $ that implies~\cref{eq:ass24} if $u_{f} = \mathcal{O}(m^{-1} n^{-1} u)$.
	
	Furthermore, we shall assume that~\cref{alg:RRRCholeskyQR} in step 2 is using strong RRQR subroutine that satisfies~\cref{eq:ass25,eq:ass26}. This can be achieved for instance  with the strong rank-revealing QR method from~\cite{gu1996efficient} with unit roundoff similar to the one required by the Givens QR i.e. $u_{f}=O(k^{-1}n^{-\frac{1}{2}}u)$. The method in~\cite{gu1996efficient}  contains an extra parameter $f$ that in our case should be taken as, say, $1.5$. Then the RRQR subroutine in step 2 will take a negligible amount $\mathcal{O}(kn^2\log n)$ of flops, while satisfying~\cref{eq:ass25,eq:ass26}.
	
	First, it is shown in~\cref{thm:RRRCholeskyQR1} that RRRCholeskyQR permutes columns so that the condition number of $\bX_{(1:r)}$ is bounded by $F(n,r) \tau^{-1}$. 
	
	\begin{theorem} \label{thm:RRRCholeskyQR1}
		Let $\bX$ have normalized columns. Consider~\cref{alg:RRRCholeskyQR} using the strong rank-revealing QR algorithm and $\tau \geq 4 n^{\frac{3}{2}} r  u $.  Let $\bTheta$ be an $\varepsilon$-embedding for $\bX_{(1:r)}$ with $\varepsilon \leq \frac{1}{2}$. Under~\cref{ass:ass200} possibly excluding~\cref{eq:ass23,eq:ass24}, we have
		\begin{equation} \label{eq:RRRCholeskyQR1}
		\cond(\bX_{(1:r)}) \leq 10 n^{\frac{3}{2}} r {\tau}^{-1}.
		\end{equation}
		\begin{proof}	
			First we notice that by~\cref{eq:ass25}, 
			\begin{equation} \label{eq:RRRCholeskyQR11}	
			\sigma_{min}(\bhR_{(1:r)}) \geq 4^{-1} n^{-1} r^{-1} \|\bhR_{(r:n)} \|_2  \geq 4^{-1} n^{-\frac{3}{2}} r^{-1} \tau \|\bhR \|_2.
			\end{equation}	
			Thus, we deduce that
			\begin{equation} \label{eq:RRRCholeskyQR112}	
			\mathrm{cond}(\bhR_{(1:r)}) \leq 4 n^{\frac{3}{2}}r \tau^{-1}\leq u^{-1}.
			\end{equation}
			
			Furthermore, by~\cref{eq:ass22} we have 
			\begin{equation} \label{eq:RRRCholeskyQR12}	
			\|\bhP_{(1:r)}\|_2 \leq \|\bhS\|_2 \|\bhR_{(1:r)}\|_2+ \|\bE^{*}_2\|_2 \leq 1.01 \|\bhR_{(1:r)}\|_2,
			\end{equation}
			and 
			\begin{equation} \label{eq:RRRCholeskyQR13}	
			\begin{split}
			\sigma_{min}(\bhP_{(1:r)}) &\geq \sigma_{min}(\bhS) \sigma_{min}(\bhR_{(1:r)}) - \|\bE^{*}_2\|_2 \geq 0.99 \sigma_{min}(\bhR_{(1:r)}) - 0.01u  \| \bP_{(1:r)} \|_2 \\ 
			& \geq  0.99 \sigma_{min}(\bhR_{(1:r)}) - 0.011 u  \| \bR_{(1:r)} \|_2 \\
			& \geq  \sigma_{min}(\bhR_{(1:r)}) (0.99 - 0.011 u \mathrm{cond}(\bhR_{(1:r)}) \geq  0.97  \sigma_{min}(\bhR_{(1:r)}).
			\end{split}
			\end{equation}
			Consequently,
			\begin{equation} \label{eq:RRRCholeskyQR14}	
			\mathrm{cond}(\bhP_{(1:r)}) \leq 4.2 n^{\frac{3}{2}} r {\tau}^{-1} \leq 1.05 u^{-1}. 
			\end{equation}			
			Next, by~\cref{eq:ass21} we get 	
			\begin{equation} \label{eq:RRRCholeskyQR15}			
			\| \bTheta \bX_{(1:r)} \|_2 \leq \| \bhP_{(1:r)} \|_2 + \|\bE^{*}_1 \|_2 \leq 1.01 \| \bhP_{(1:r)} \|_2, 
			\end{equation}
			which due to the $\varepsilon$-embedding property of $\bTheta$ implies that 
			\begin{equation} \label{eq:RRRCholeskyQR16}			
			\| \bX_{(1:r)} \|_2 \leq 1.5 \| \bhP_{(1:r)} \|_2.
			\end{equation}			
			We also have by~\cref{eq:ass21,eq:RRRCholeskyQR14},
			\begin{equation} \label{eq:RRRCholeskyQR17}			
			\begin{split}
			\sigma_{min}(\bTheta \bX_{(1:r)}) &\geq \sigma_{min}(\bhP_{(1:r)}) - \|\bE^{*}_1\|_2 \geq  \sigma_{min}(\bhP_{(1:r)}) -0.01u {\| \bX_{(1:r)} \|_2}  \\ 
			& \geq  \sigma_{min}(\bhP_{(1:r)}) -0.015u \|\bhP_{(1:r)}\|_2 
			\geq  0.92 \sigma_{min}(\bhP_{(1:r)}).
			\end{split}
			\end{equation}			
			Consequently, by the $\varepsilon$-embedding property of $\bTheta$ and~\cref{eq:RRRCholeskyQR16,eq:RRRCholeskyQR17},
			\begin{equation*}
			\mathrm{cond}(\bX_{(1:r)}) \leq \sqrt{\frac{1+\varepsilon}{1-\varepsilon}}\mathrm{cond}(\bTheta\bX_{(1:r)}) \leq 1.91 \mathrm{cond}(\bhP_{(1:r)}) \leq  10 n^{\frac{3}{2}} r {\tau} ^{-1},
			\end{equation*}
			that finishes the proof.
		\end{proof}	
	\end{theorem}
	
	\Cref{thm:RRRCholeskyQR1} implies stability of the computation of the Q factor by forward substitution in step 5. With this result we are ready to establish the overall stability guarantee. 
	
	\begin{theorem} \label{thm:RRRCholeskyQR2}
		Let $\bX$ have normalized columns. Consider~\cref{alg:RRRCholeskyQR} using the strong rank-revealing QR algorithm and $1000 r^\frac{5}{2} n^{\frac{3}{2}} u \leq  {\tau} \leq 1$. Let $\bTheta$ be an $(\varepsilon, \delta, n)$-OSE with $\varepsilon \leq \frac{1}{2}$. Under~\cref{ass:ass200}, we have with probability at least $1-\delta$,
		\begin{subequations} \label{eq:RRRCholeskyQR2}
			\begin{align}	
			\|\bX - \bhQ \bhR\|_\Frob &\leq 2 {\tau}   , \label{eq:main4} \\
			(1+\varepsilon)^{-\frac{1}{2}} - 0.016  \leq \sigma_{min}(\bhQ)&\leq \sigma_{max}(\bhQ)  \leq (1-\varepsilon)^{-\frac{1}{2}}+ 0.016, \label{eq:main5}
			\end{align}
			In addition, it holds that
			\begin{equation} \label{eq:main6}
			\|\bhS - \bTheta \bhQ \|_\Frob \leq 61 r^\frac{5}{2} n^{\frac{3}{2}} \frac{u}{{\tau}} \leq 0.061.
			\end{equation}
		\end{subequations}
		Furthermore, the stability guarantees~\cref{eq:RRRCholeskyQR2} hold  with probability at least $1-\delta$, if $\bTheta$ is an $(\varepsilon, \binom{n}{r+1}^{-1} \delta, r+1)$-OSE and not necessarily an $(\varepsilon, \delta, n)$-OSE. 
		\begin{proof}
			Similarly as in~\cref{thm:rrCholeskyQR}, we here shall assume that $\bTheta$ is an $\varepsilon$-embedding for all subspaces spanned by $r+1$ columns of $\bX$. As is argued in the proof of~\cref{thm:rrCholeskyQR} this condition is satisfied with probability at least $1-\delta$. Then by~\cref{eq:ass21,eq:ass22,eq:ass24,eq:ass26} we have	
			\begin{equation} \label{eq:RRRCholeskyQR20}
			\begin{split}
			\| \bX_{(r+1:n)}  - \bhQ \bhR_{(r+1:n)} \|_\Frob &=\| \bX_{(r+1:n)}  - \bX_{(1:r)} \bhR_{(1:r)}^{-1} \bhR_{(r+1:n)}\|_\Frob + \| \bE^{*}_3 \bhR_{(1:r)}^{-1} \bhR_{(r+1:n)} \|_\Frob \\&\leq (1-\varepsilon)^{-\frac{1}{2}} \|\bTheta(\bX_{(r+1:n)}  - \bX_{(1:r)} \bhR_{(1:r)}^{-1} \bhR_{(r+1:n)}) \|_\Frob + \| \bE^{*}_3\|_\Frob \|\bhR_{(1:r)}^{-1} \bhR_{(r+1:n)} \|_2\\ &\leq (1-\varepsilon)^{-\frac{1}{2}}( \|\bE_6 \|_\Frob +  \|\bE_7 \|_\Frob + \|\bE_8 \|_\Frob+ \|\bE_{9} \|_\Frob ) + 2n^\frac{1}{2} r^\frac{1}{2}\| \bE^{*}_3 \|_\Frob,
			\end{split}
			\end{equation}
			where
			\begin{align*}
			\|\bE_6 \|_\Frob &= \| \bTheta \bX_{(r+1:n)} - \bhP_{(r+1:n)} \|_\Frob \leq 0.1 u \|  \bX_{(r+1:n)} \|_2  \\
			\|\bE_7 \|_\Frob &= \|\bhP_{(r+1:n)} - \bhS \bhR_{(r+1:n)}\|_\Frob \leq 1.01 {\tau} \|\bhR\|_2 \leq 1.02 (1+\varepsilon)^{\frac{1}{2}}{\tau} \|\bX\|_\Frob \\
			\|\bE_{8} \|_\Frob &= \|(\bhS \bhR_{(1:r)}-\bhP_{(1:r)}) \bR_{(1:r)}^{-1} \bR_{(r+1:n)}\|_\Frob \leq \|\bhS \bhR_{(1:r)}-\bhP_{(1:r)}\|_\Frob \|\bR_{(1:r)}^{-1} \bR_{(r+1:n)}\|_2  \\ &\leq 2 n^\frac{1}{2} r^\frac{1}{2} \|\bE^{*}_2 \|_\Frob \leq 0.02 n^\frac{1}{2} r^\frac{1}{2} u \|\bhP_{(1:r)} \|_2 \leq 0.03 n u \|\bX_{(1:r)} \|_2
			\\
			\|\bE_9 \|_\Frob & = \| ( \bhP_{(1:r)} - \bTheta \bX_{(1:r)} ) \bhR_{(1:r)}^{-1} \bhR_{(r+1:n)}\|_2\leq \| \bTheta \bX_{(1:r)} - \bhP_{(1:r)}  \|_\Frob \| \bhR_{(1:r)}^{-1} \bhR_{(r+1:n)}\|_2 \\& \leq 2 n^\frac{1}{2} r^\frac{1}{2} \| \bE^{*}_1 \|_\Frob \leq 0.02 n u \|  \bX_{(1:r)} \|_2.
			\end{align*}		
			Consequently,
			\begin{equation} \label{eq:RRRCholeskyQR21}
			\| \bX_{(r+1:n)}  - \bhQ \bhR_{(r+1:n)} \|_\Frob \leq \sqrt{2} ( un + 1.02 \sqrt{\textstyle\frac{3}{2}}{\tau} )\|\bX\|_\Frob.
			\end{equation}
			Furthermore, from~\cref{thm:RRRCholeskyQR1} it follows that 
			\begin{equation} \label{eq:RRRCholeskyQR31}
			\cond(\bX_{(1:r)}) \leq 10 n^{\textstyle \frac{3}{2}} r {\tau}^{-1}.
			\end{equation}
			By looking at $\bhQ$ and $\bhR_{(1:r)}$ as a RCholeskyQR factorization of $\bX_{(1:r)}$, according to~\cref{thm:main}, we have 
			\begin{subequations}
				\begin{align}	
				\|\bX_{(1:r)} - \bhQ \bhR_{(1:r)}\|_\Frob &\leq 2.1 u r \|\bX_{(1:r)} \|_\Frob \label{eq:RRRmain1} \\
				(1+\varepsilon)^{-\frac{1}{2}} - 4 u r^\frac{3}{2}  \mathrm{cond}(\bX_{(1:r)}) \leq \sigma_{min}(\bhQ)&\leq \sigma_{max}(\bhQ)  \leq (1-\varepsilon)^{-\frac{1}{2}}+ 4 u r^\frac{3}{2}  \mathrm{cond}(\bX_{(1:r)}) \label{eq:RRRmain2} \\
				\label{eq:RRRmain3}
				\|\bhS - \bTheta \bhQ \|_\Frob &\leq 6.1 u r^\frac{3}{2}  \mathrm{cond}(\bX_{(1:r)}).
				\end{align}
			\end{subequations}
			
			By combing~\cref{eq:RRRmain1} with~\cref{eq:RRRCholeskyQR21} we obtain~\cref{eq:main4}. By combining~\cref{eq:RRRmain2} with~\cref{eq:RRRCholeskyQR31} we obtain~\cref{eq:main5}. Finally, by combing~\cref{eq:RRRmain3} with~\cref{eq:RRRCholeskyQR31} we obtain~\cref{eq:main6} and finish the proof. 	
		\end{proof}

	\end{theorem}
	
	\section{Numerical experiments} \label{numexperiments}
	
	In this section, the proposed RCholeskyQR and RRRCholeskyQR factorizations are verified on numerical examples. They are compared in terms of stability with the standard CholeskyQR2 and Householder QR factorizations, as well as with the shifted CholeskyQR2 and shifted CholeskyQR3 from~\cite{fukaya2020shifted}, the classical block Gram-Schmidt algorithm with reorthogonalization (BCGS2), the modified block Gram-Schmidt algorithm (BMGS) and the block RGS algorithm from~\cite {balabanov2021randomized}. We also characterize the potential speedups of RCholeskyQR2 and RRRCholeskyQR2 methods over the shifted CholeskyQR3 and Householder QR.

	\subsection{Direct QR factorization}
	The algorithms were tested on series $\bX^{(1)}, \bX^{(2)}, \hdots, \bX^{(j)}$ of matrices of various types and sizes.  For the sake of completeness, we have considered two scenarios: one that concerns obtaining a well-conditioned Q factor with RCholeskyQR, RRRCholeskyQR, CholeskyQR2, shifted CholeskyQR2, or Householder QR, and one that concerns obtaining an orthonormal Q factor with RCholeskyQR2, RRRCholeskyQR2, shifted CholeskyQR3, or Householder QR. On the plots, these methods are denoted by RCholQR, RRRCholQR, CholQR2, sCholQR2, HH and RCholQR2, RRRCholQR2, sCholQR3, HH, respectively. In the first case, stability was characterized by $\cond(\bQ)$, and in the second case, by the classical measure $\Delta = \|\bQ^{\mathrm{T}}\bQ - \bI\|_2$. Since in randomized algorithms the SRHT and Gaussian matrices gave very similar results, here we present the results for SRHT only.

	To ensure a fair comparison, before executing the shifted CholeskyQR algorithms, we normalized the columns of $\bX^{(i)}$. Furthermore, we tested several variants of shifts, and then chose those that gave the greatest stability. Namely, the first orthogonalization was performed considering the shift $s$ as the maximum between $s_0$ and the smallest power of $10$, such that $\bX^\mathrm{T}\bX+ s \bI$ is numerically positive-definite. The $s_0$ parameter was chosen to be either zero, or the recommended value  from~\cite{fukaya2020shifted}: $11u (mn+n(n+1))\| \bX\|^2_2$, or empirically chosen value: $u \sqrt{n} \| \bX\|^2_\Frob$. In addition, the second orthogonalization was performed either with zero shift or the smallest power of $10$, so that $\bQ^\mathrm{T}\bQ$ is numerically positive definite. In the shifted CholeskyQR3, the third orthogonalization was performed with zero shift.
	
	In the RCholeskyQR algorithm, we took the $\mathtt{QR}$ subroutine in step 2 as the  Householder QR.  In RRRCholeskyQR, we took $\mathtt{RRQR}$ as a strong rank-revealing QR from~\cite{gu1996efficient} with parameter $f=1.5$. Furthermore, we considered an (implicit) orthonormalization of the columns of $\bX^{(i)}$ as explained in~\cref{rrrandCholeskyQR}.  The truncation parameter $\tau$ was chosen to be of order of $10^{-15}$ for experiments in float64 format and of order of $10^{-7}$ for experiments in float32. It was revealed that, depending on the experiment, changing this parameter by a small factor could slightly improve the results. Even though the improvements were not that significant, for a fair comparison, we decided to present the results corresponding to the best tested values of $\tau$.

\begin{figure}[!h]
		\centering
		\begin{subfigure}{.35\textwidth}
			\centering  
			\includegraphics[width=\textwidth]{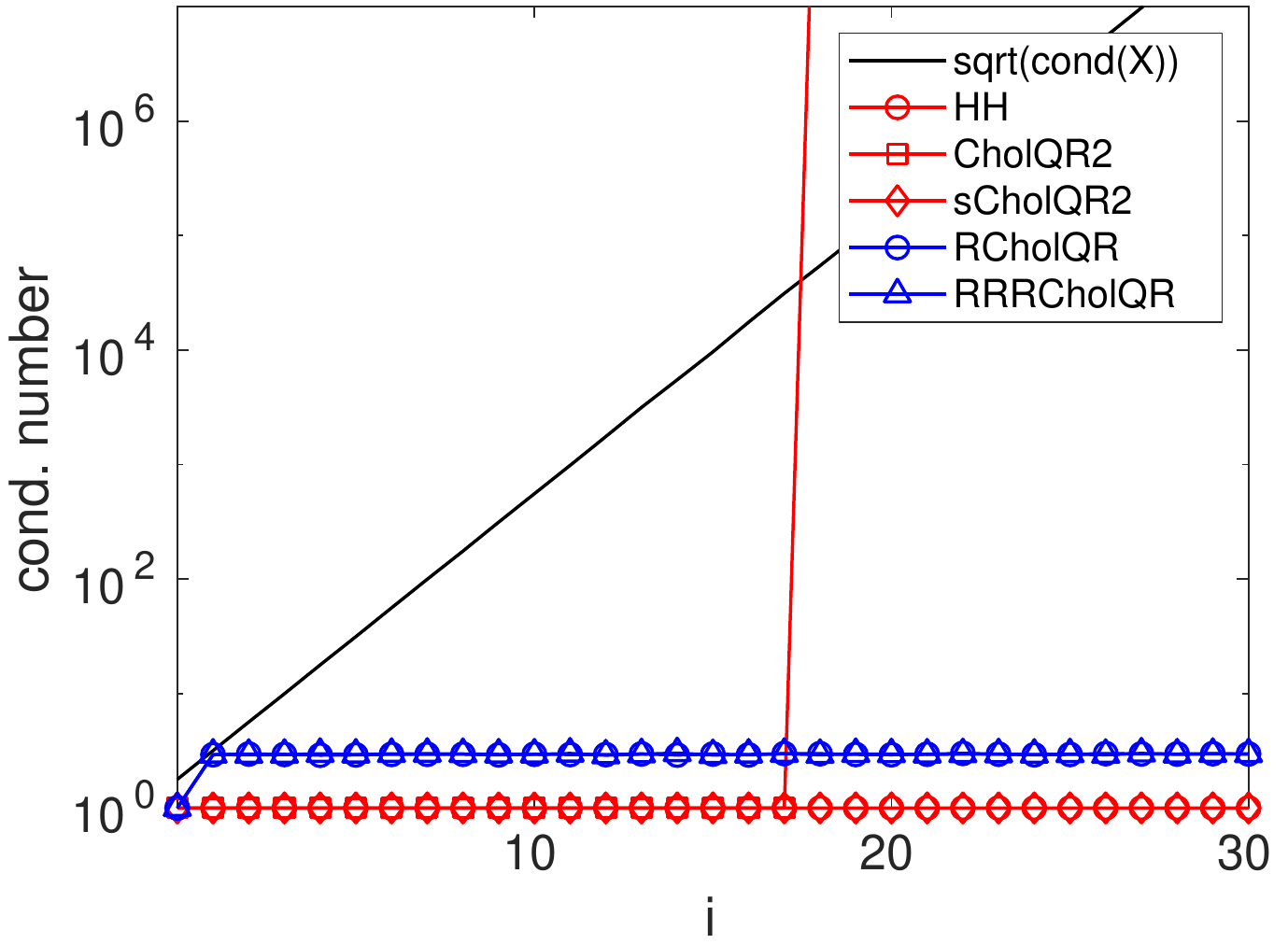}
			\caption{\small Cond. number of $\bQ$.}
			\label{fig:Ex1_1a}
		\end{subfigure} \hspace{.03\textwidth}
		\begin{subfigure}{.35\textwidth}
			\centering
			\includegraphics[width=\textwidth]{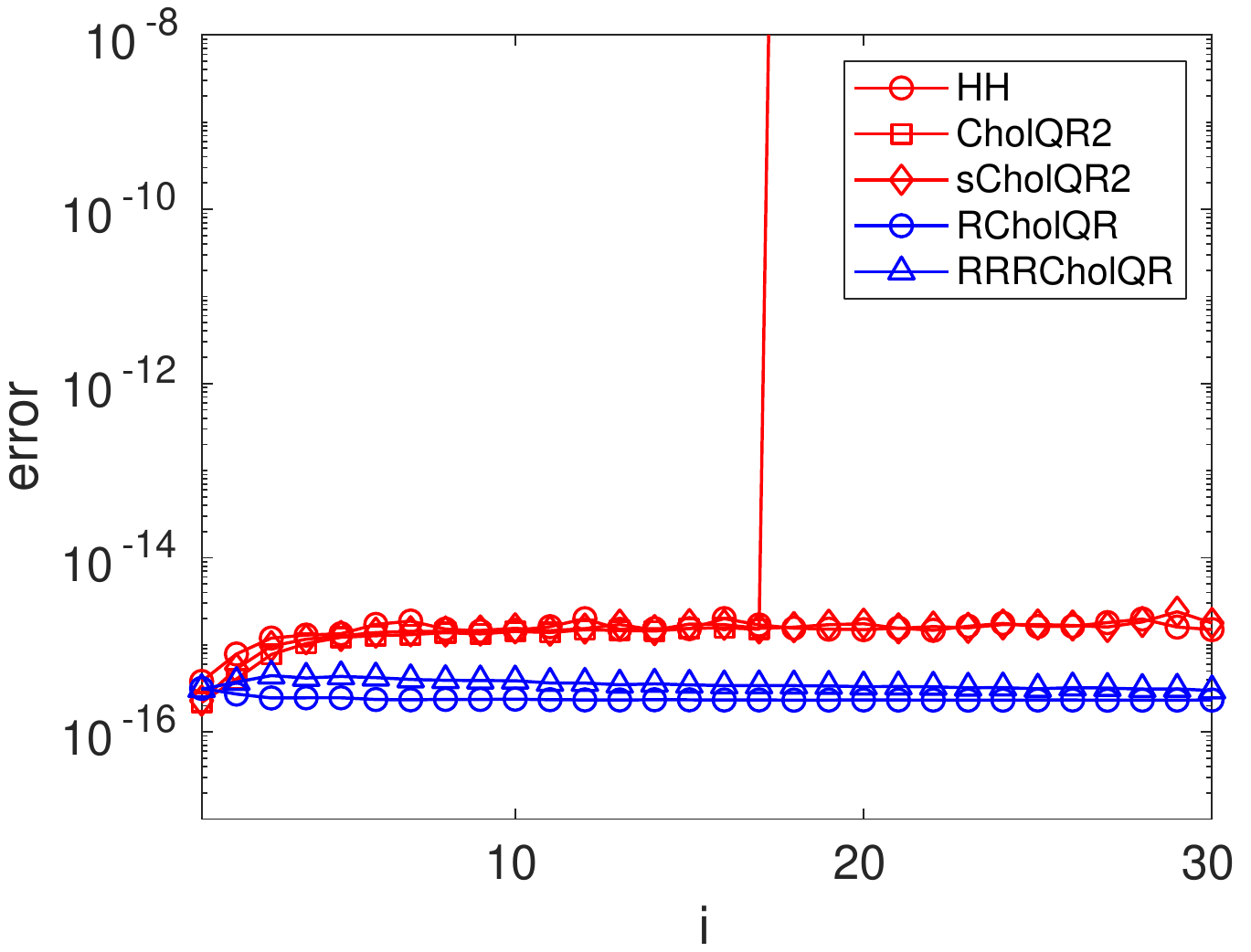}
			\caption{\small Max relative column-wise error.}
			\label{fig:Ex1_1b}
		\end{subfigure} \\
		\begin{subfigure}{.35\textwidth}
			\centering  
			\includegraphics[width=\textwidth]{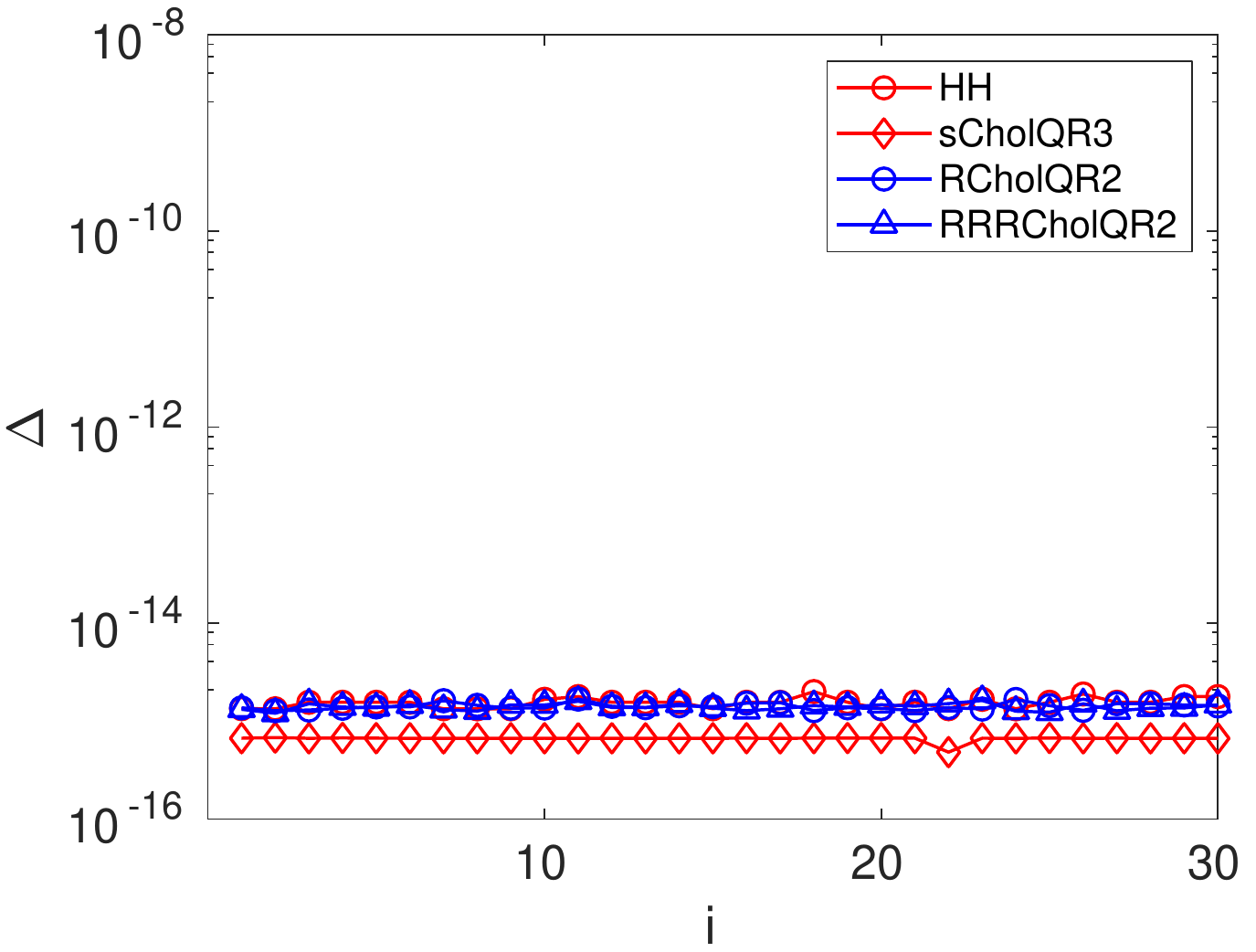}
			\caption{\small Stability measure  $\Delta = \|\bQ^{\mathrm{T}}\bQ - \bI\|_2$.}
			\label{fig:Ex1_1c}
		\end{subfigure} \hspace{.03\textwidth}
		\begin{subfigure}{.35\textwidth}
			\centering
			\includegraphics[width=\textwidth]{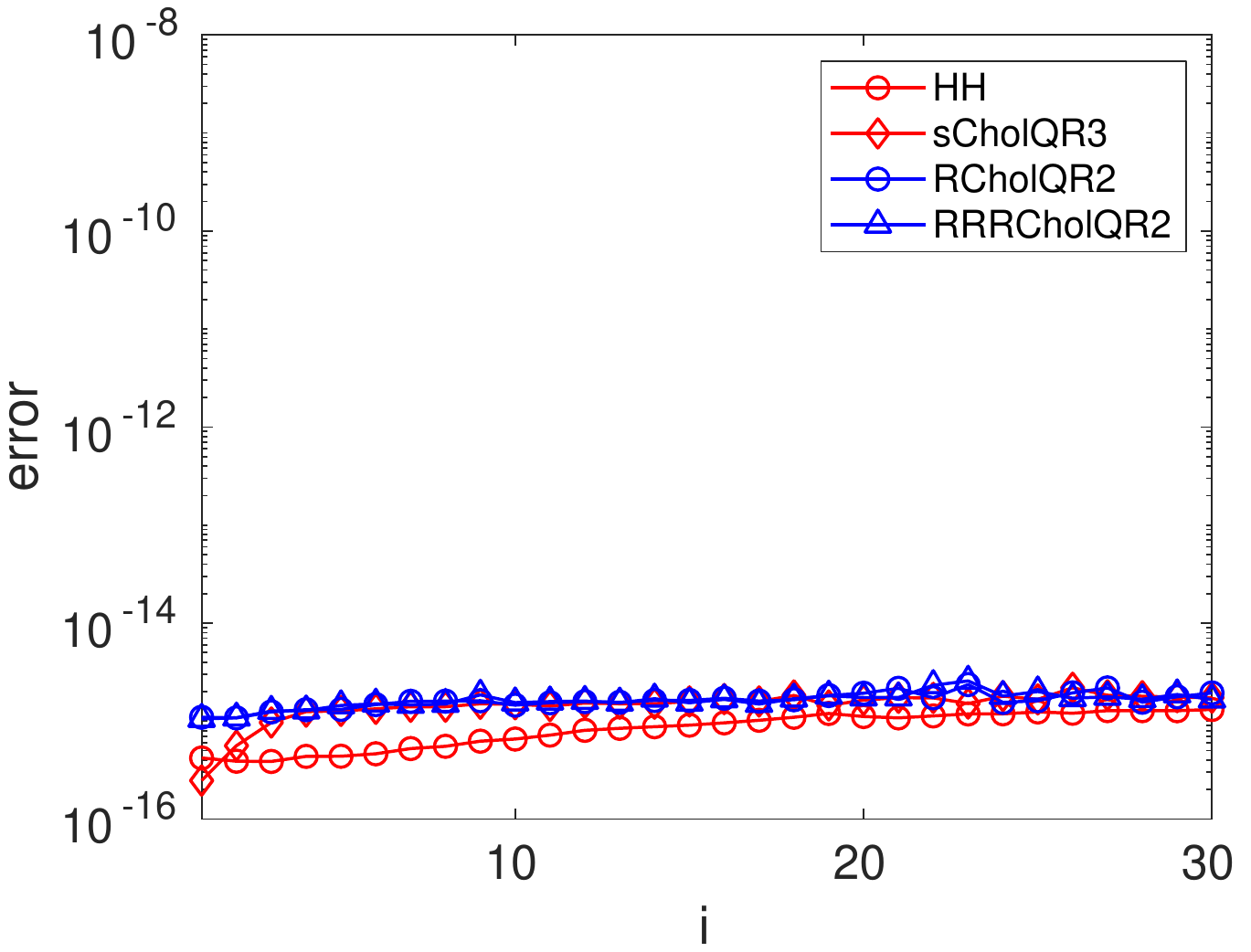}
			\caption{\small Max relative column-wise error.}
			\label{fig:Ex1_1d}
		\end{subfigure}	
		\caption{Stability characterization of QR factorizations of $\bX^{(i)}$ of the form $\bU \bSigma \bV^\mathrm{T}$, where $\bU$ and $\bV $ are Gaussian matrices orthonormalized by Householder QR, and $\bSigma = \mathrm{diag}([1,\sigma^{\frac {1}{ n-1}}, \hdots, \sigma^{\frac{n-2}{n-1}},\sigma])$, with parameter $\sigma=\sigma(i)$ ranging from $10^ {-15}$ to $1$.}
		\label{fig:Ex1_1}
	\end{figure}

	In the first test, the matrices $\bX^{(i)}$ were taken of the form $\bX^{(i)} = \bU \bSigma \bV^\mathrm{T}$, where $\bU$ and $\bV $ are $m\times n$ and $n \times n$ random Gaussian matrices orthonormalized by the Householder QR, and $\bSigma = \mathrm{diag}([1,\sigma^{\frac {1}{ n-1}}, \hdots, \sigma^{\frac{n-2}{n-1}},\sigma])$, with parameter $\sigma=\sigma(i)$ ranging from $10^ {-15}$ to $1$, controlling the condition number of $\bX^{(i)}$, as in~\cite{fukaya2020shifted}. The dimensions $m$ and $n$ were chosen as $10^6$ and $300$ respectively, and the sketching size $k$ as $2n=600$. In this experiment, all operations were performed in float64 format with unit rounding $\approx 10^{-16}$.  In RRRCholeskyQR we chose the $\tau$ parameter to be $4\times 10^{-15}$.
	We first turned to computing the QR factorization with a well-conditioned Q factor. From~\cref{fig:Ex1_1a,fig:Ex1_1b} one can see that RCholeskyQR and RRRCholeskyQR were very stable for all $\bX^{(i)}$, as were the shifted CholeskyQR2 and Householder QR. The standard CholeskyQR2, on the other hand, failed when $\mathrm{cond}(\bX^{(i)})$ got larger than about $10^{8} \approx u^{\frac{1}{2}}$, which is in good agreement with the theory. In~\cref{fig:Ex1_1c,fig:Ex1_1d} we depict the accuracy of randomized and deterministic QR factorizations in terms of orthogonality of the Q factor. As in the previous experiment, RCholeskyQR2, RRRCholeskyQR2, shifted CholeskyQR3 and Householder QR provided near perfect stability. 
	
\begin{figure}[!h] \label{fig:Ex2}
		\centering
		\begin{subfigure}{.35\textwidth}
			\centering  
			\includegraphics[width=\textwidth]{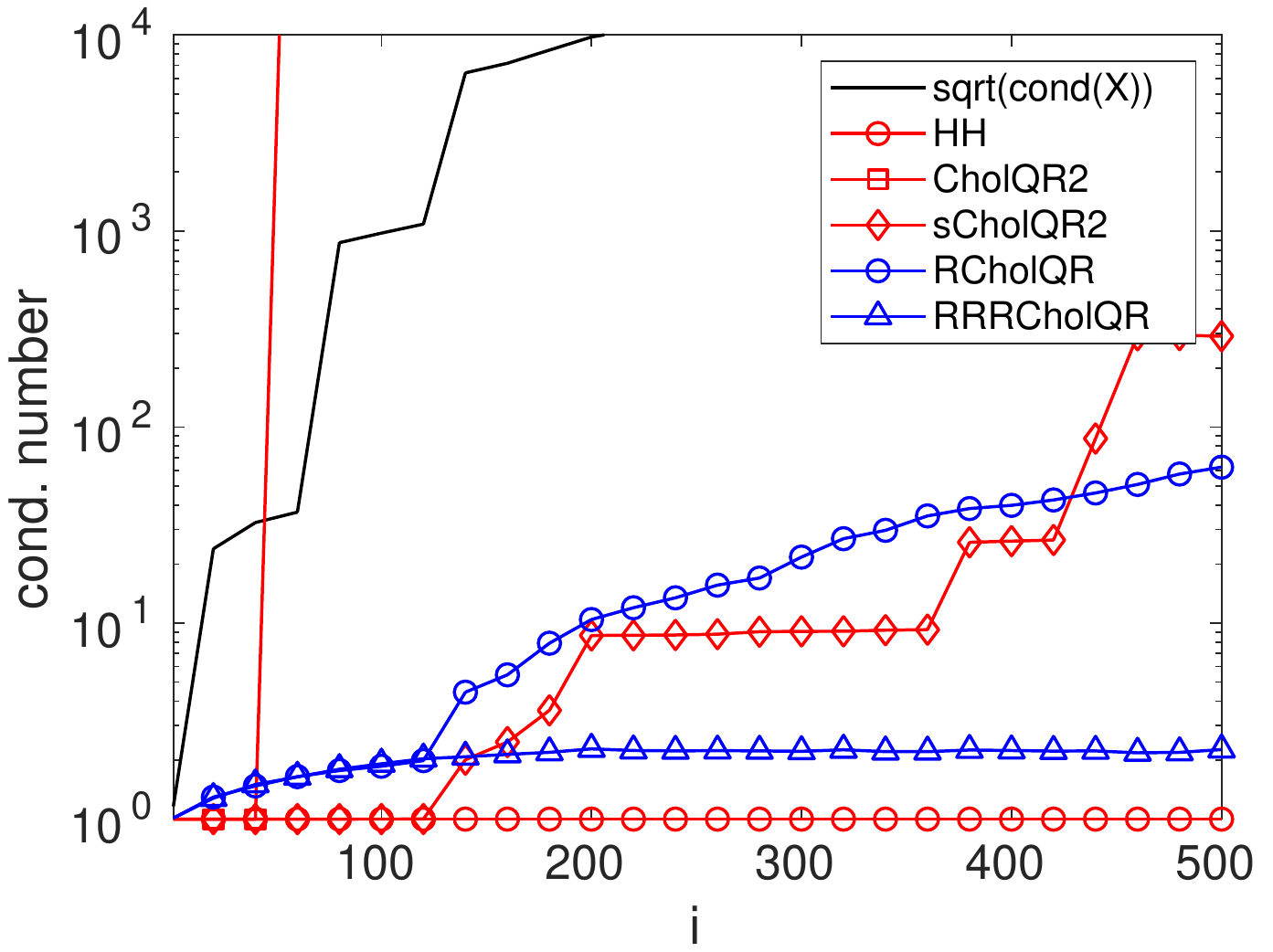}
			\caption{\small Cond. number of $\bQ$.}
			\label{fig:Ex2_1a}
		\end{subfigure} \hspace{.03\textwidth}
		\begin{subfigure}{.35\textwidth}
			\centering
			\includegraphics[width=\textwidth]{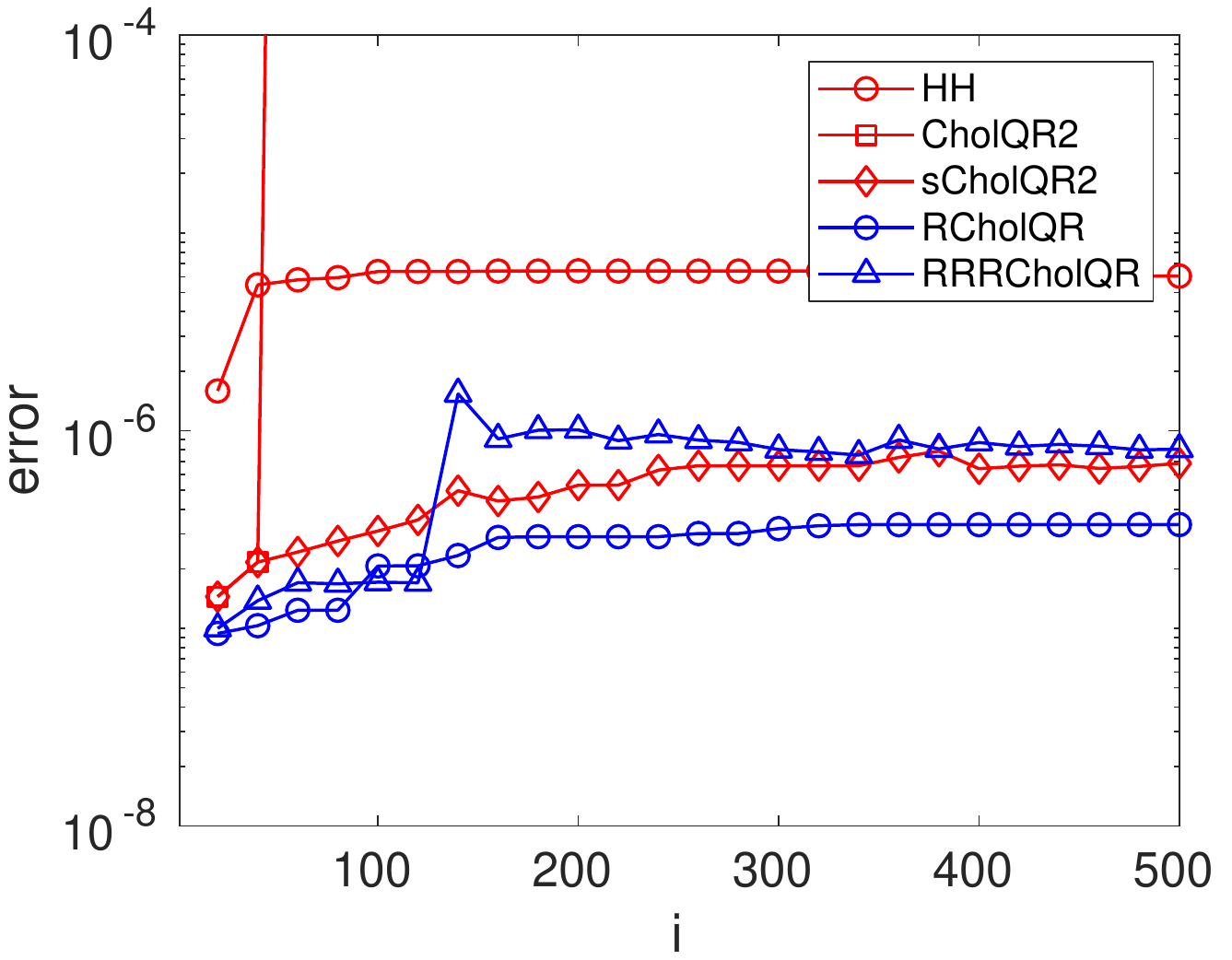}
			\caption{\small Max relative column-wise error.}
			\label{fig:Ex2_1b}
		\end{subfigure} \\
		\begin{subfigure}{.35\textwidth}
			\centering  
			\includegraphics[width=\textwidth]{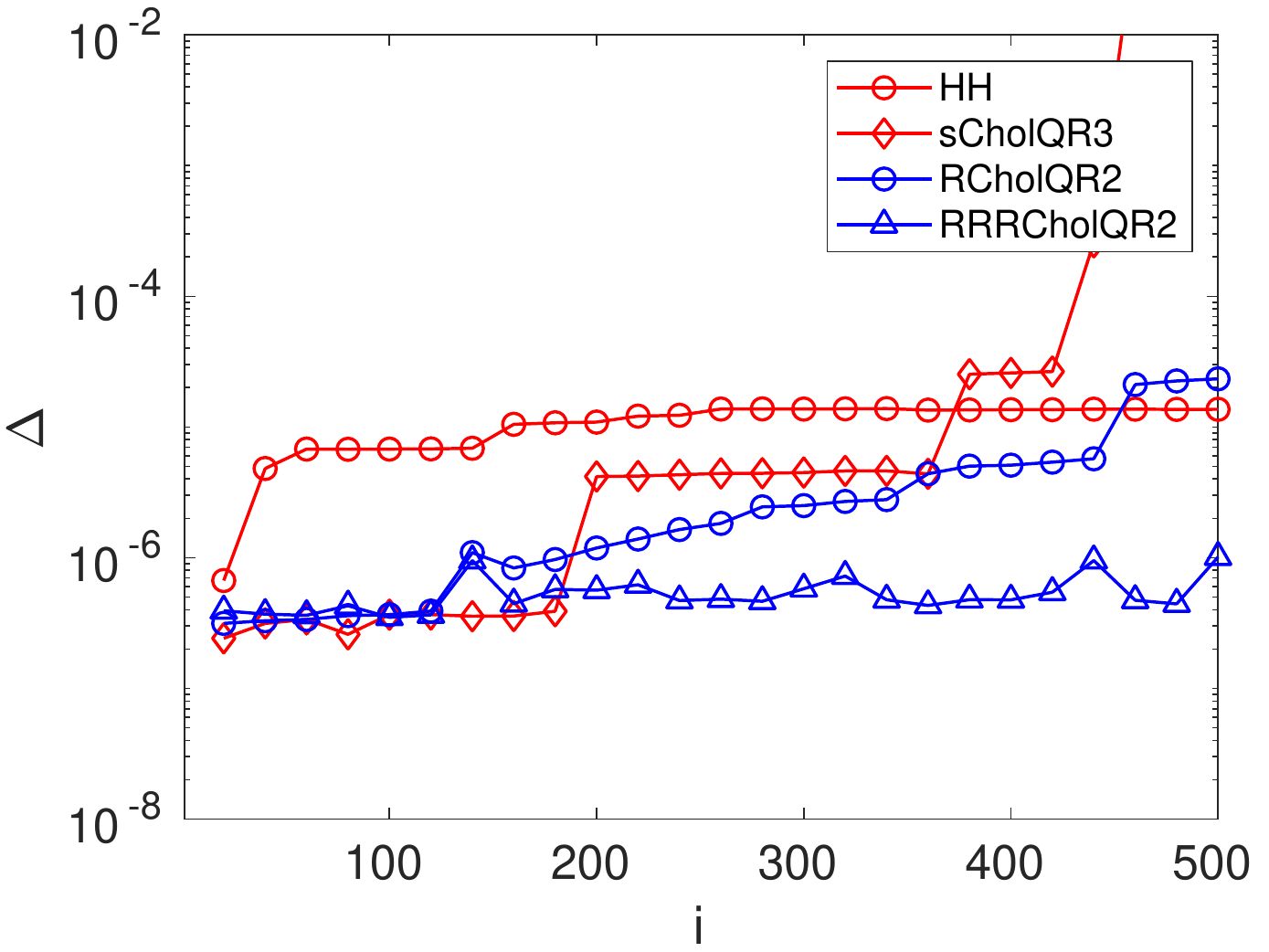}
			\caption{\small Stability measure  $\Delta = \|\bQ^{\mathrm{T}}\bQ - \bI\|_2$.}
			\label{fig:Ex2_1c}
		\end{subfigure} \hspace{.03\textwidth}
		\begin{subfigure}{.35\textwidth}
			\centering
			\includegraphics[width=\textwidth]{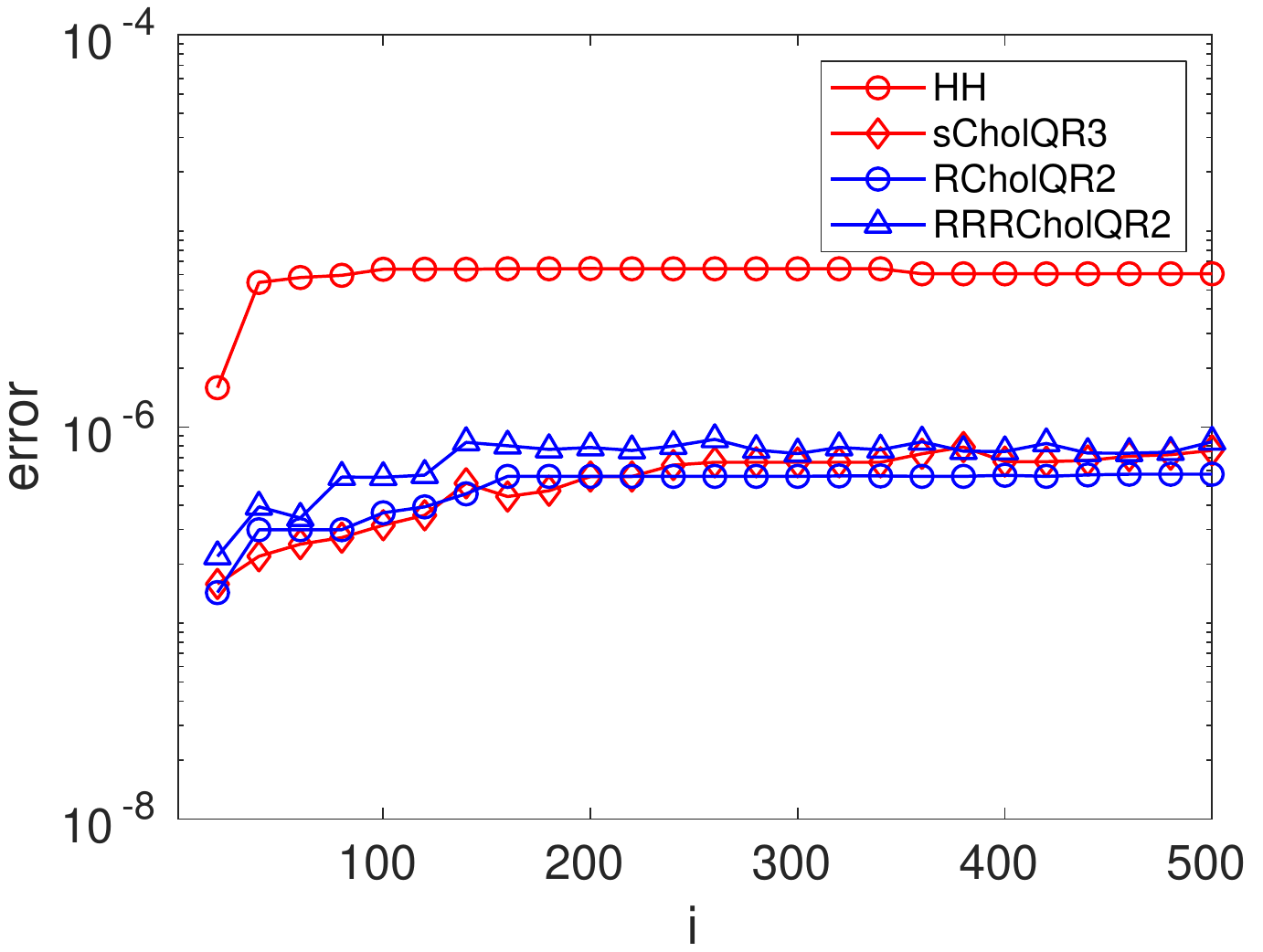}
			\caption{\small Max relative column-wise error.}
			\label{fig:Ex2_1d}
		\end{subfigure}	
		\caption{{Stability characterization of QR factorizations of $\bX^{(i)} = \bW_{(1:i)}$, where $\bW$ is a discretization of $f(x,y)$ on an uniform $10^6\times 500$ grid. }}
		\label{fig:Ex2_1}
	\end{figure}

	In the second test case, we constructed the $\bX^{(i)}$ matrices as in~\cite{balabanov2021randomizedGS}. In particular, we took an uniform unit grid of size $m \times n$ with $m = 10^6$ and $n=500$, and constructed an $m \times n$ matrix $\bW$ with entries equal to 
	\begin{equation*}
	f(x,y) = \frac{\sin\left (10(y+x) \right)}{\cos \left(100(y-x) \right)+1.1}
	\end{equation*}
	evaluated at the corresponding grid points. Then we considered QR factoizations of matrices $\bX^{(i)} = \bW_{(1:i)}$, $1\leq i \leq n$. In randomized algorithms the matrix $\bTheta$ was taken of size $k=2n=1000$. Here we used float32 arithmetic with working roundoff $u\approx 10^{-7}$. In RRRCholeskyQR the $\tau$ parameter was taken as $2\times 10^{-7}$. In addition, we tested the benefits of using multi-precision arithmetic in randomized algorithms. It turned out that performing minor operations in the float64 format allowed to reduce the condition number of the Q factor in RRRCholeskyQR by almost an order of magnitude, while in RCholeskyQR the multi-precision framework did not provide a significant advantage. Therefore, the following results for RRRCholeskyQR will be for the multi-precision algorithm, and the results for RCholeskyQR will be for the unique precision algorithm. Similar considerations are valid for RRRCholeskyQR2 and RCholeskyQR2. From~\cref{fig:Ex2_1}, we see that all tested methods besides the standard CholeskyQR2 showed great stability when $\bX^{(i)}$ were full-rank, which is in good agreement with the theory. 
	However, when $\bX^{(i)}$ became numerically rank-deficient at $i \geq 110$ the stability of RCholeskyQR and shifted CholeskyQR2 became deteriorated (see~\cref{fig:Ex2_1a}). Note that for RCholeskyQR the instabilities were not as high as for shifted CholeskyQR2.    In contrast, the RRRCholeskyQR algorithm provided a near-perfect stability for all $\bX^{(i)}$.  This supports the guarantee of an unconditional stability of RRRCholeskyQR from~\cref{stabRRRCholeskyQR}. Moreover, judging by the approximation errors, this algorithm turned out to be even more accurate than the Householder QR, which, we believe, is a consequence of the use of a multi-precision arithmetic framework. Now let us turn to the context of obtaining an orthonormal Q factor. According to~\cref{fig:Ex2_1c}, in this case RRRCholeskyQR2 showed greater stability than Householder QR and, in particular, provided a lower measure of stability $\|\bQ^{\mathrm{T}}\bQ - \bI\|_2$ by more than an order of magnitude. It is worth noting that RCholeskyQR2 surprisingly also showed similar or better stability as the Householder QR, while the shifted CholeskyQR3 for some $\bX^{(i)}$ provided a higher stability measure by more than two orders of magnitude. 

\begin{figure}[!h]
		\centering
		\begin{subfigure}{.35\textwidth}
			\centering  
			\includegraphics[width=\textwidth]{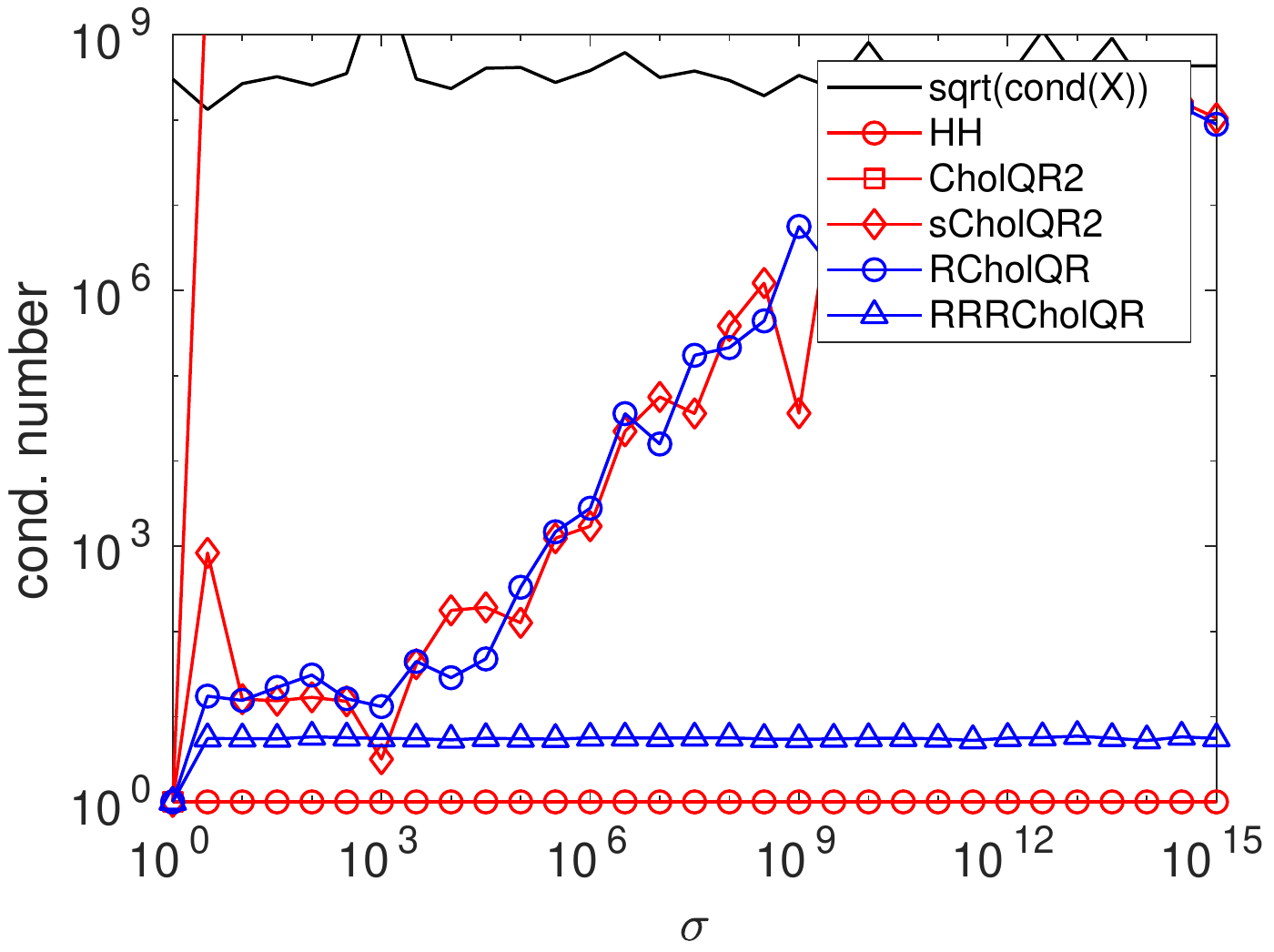}
			\caption{\small Cond. number of $\bQ$.}
			\label{fig:Ex3_1a}
		\end{subfigure} \hspace{.03\textwidth}
		\begin{subfigure}{.35\textwidth}
			\centering
			\includegraphics[width=\textwidth]{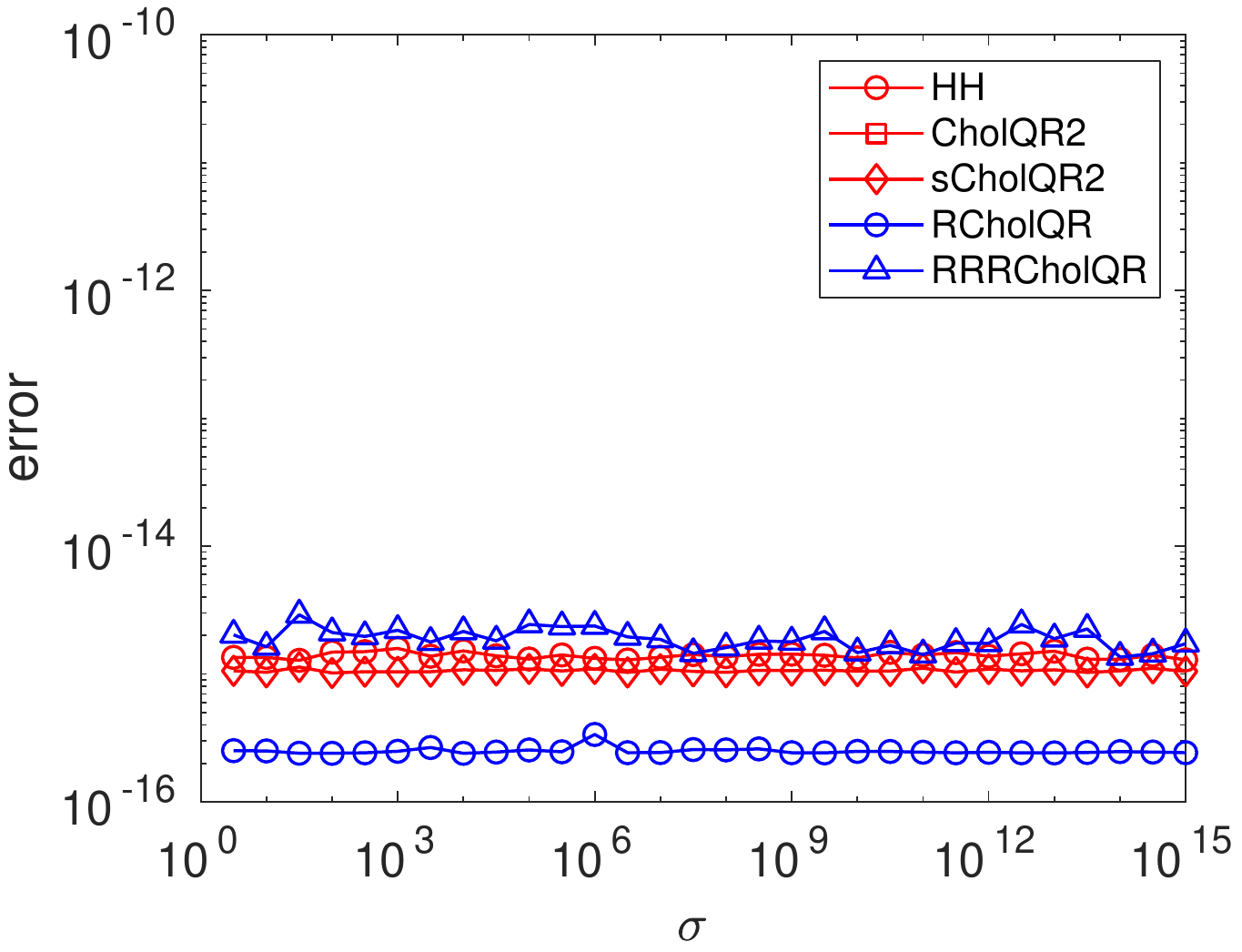}
			\caption{\small Max relative column-wise error.}
			\label{fig:Ex3_1b}
		\end{subfigure} \\
		\begin{subfigure}{.35\textwidth}
			\centering  
			\includegraphics[width=\textwidth]{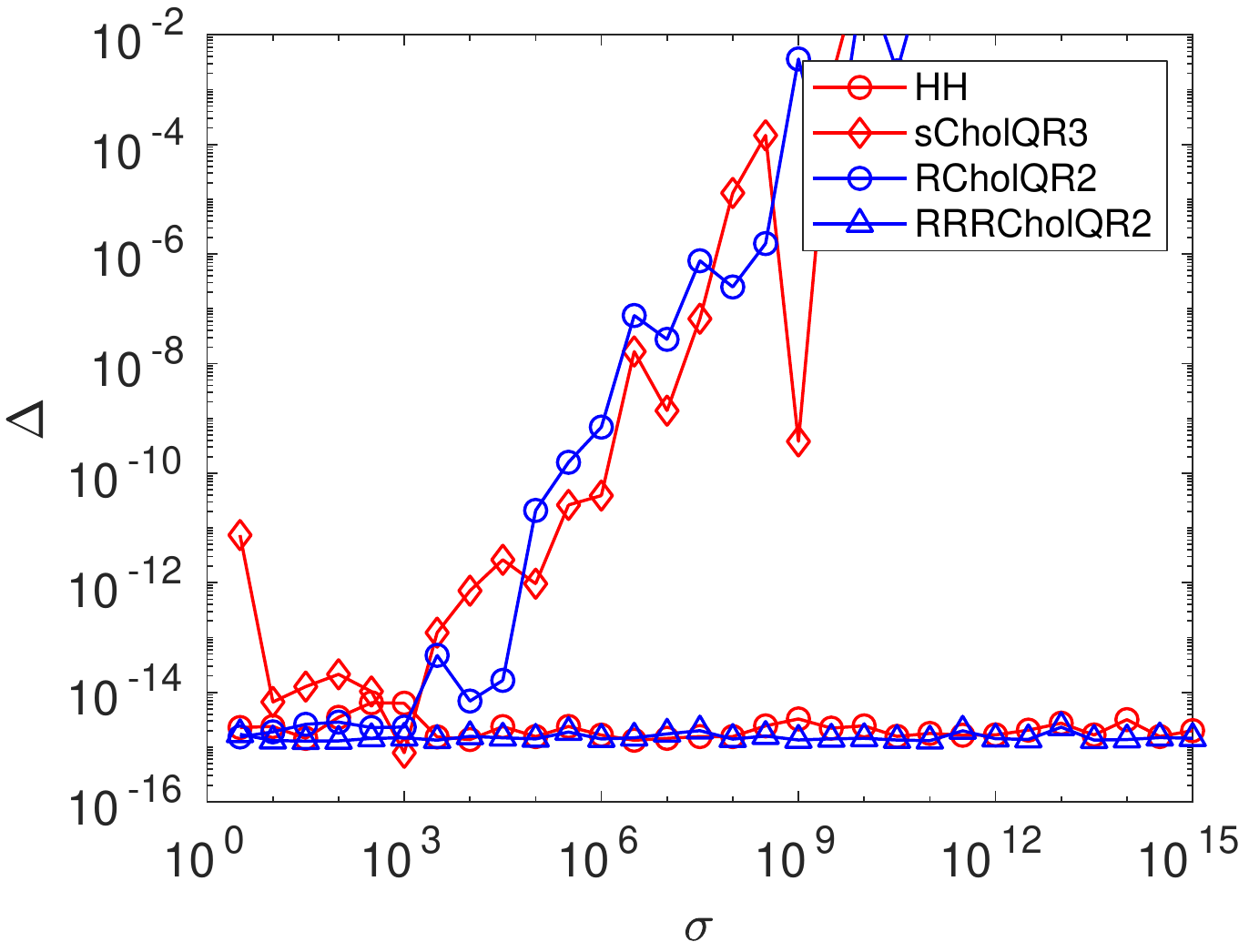}
			\caption{\small Stability measure  $\Delta = \|\bQ^{\mathrm{T}}\bQ - \bI\|_2$.}
			\label{fig:Ex3_1c}
		\end{subfigure} \hspace{.03\textwidth}
		\begin{subfigure}{.35\textwidth}
			\centering
			\includegraphics[width=\textwidth]{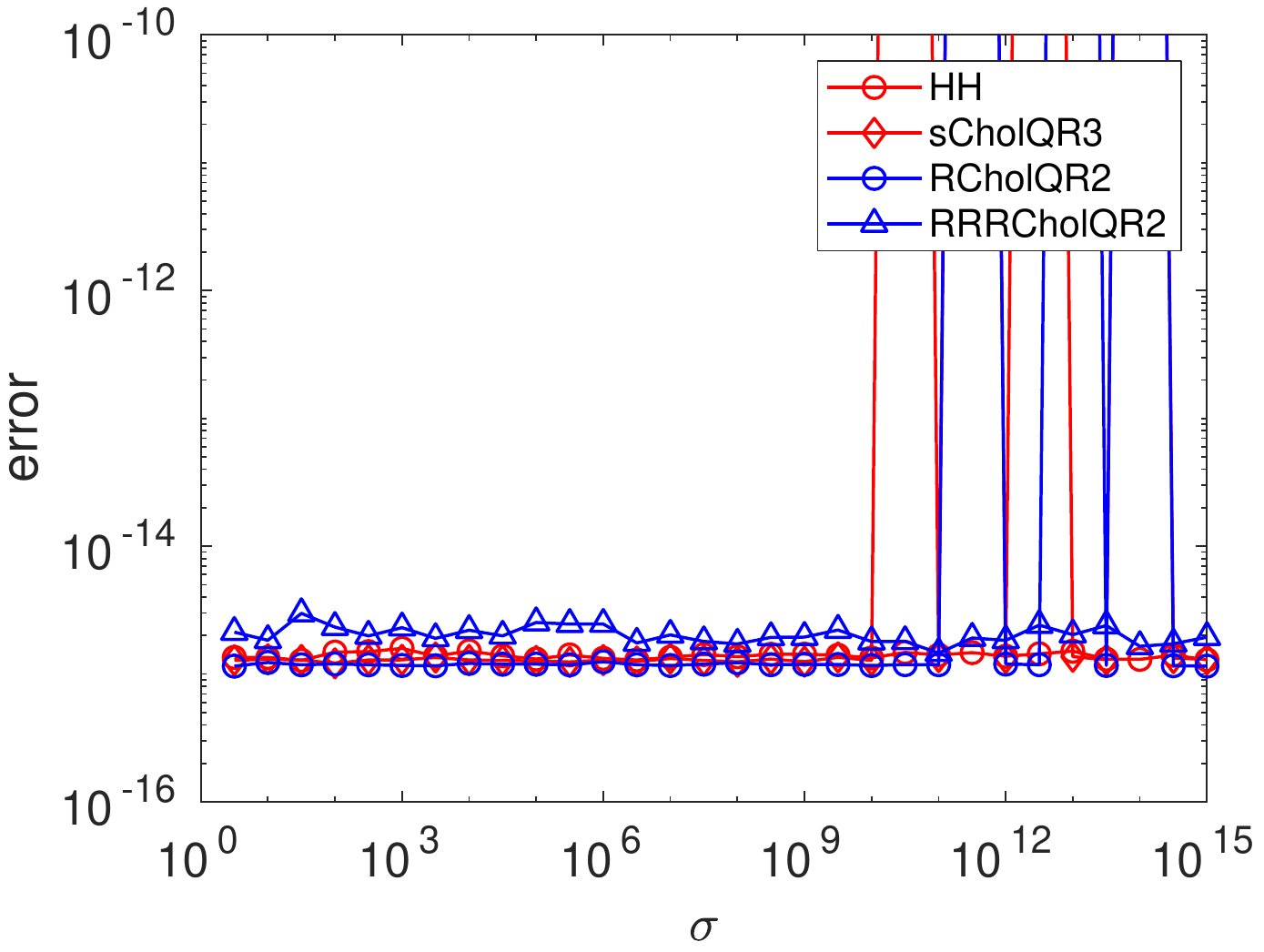}
			\caption{\small Max relative column-wise error.}
			\label{fig:Ex3_1d}
		\end{subfigure}	
		\caption{{Stability characterization of QR factorizations of rank-deficient $\bX^{(i)} = \bU \bV$ where $\bU$ is a Gaussian matrix whose first row was scaled by $\sigma$, and that was then orthonormalized by Householder QR. }}
		\label{fig:Ex2=3_1}
	\end{figure}

	To validate the unconditional stability of RRRCholeskyQR to a better extent, the third experiment involves QR factorization of rank-deficient test matrices that are particularly poorly suited to QR factorization from a numerical point of view. We have taken $\bX^{(i)}$ of the form $\bU \bV$, where the matrix $\bU$ is an $m \times n$ Gaussian matrix whose first row was scaled by a factor $\sigma(i)$ in the range from $1$ to $10^{15}$, and that was orthonormalized with Householder QR. The matrix $\bV$ is the upper triangular part of an $n \times n$ orthonormalized Gaussian matrix with modified diagonal entries to $\mathrm{diag}(\bV) = [1,10^{-15}, \hdots, 10^{-15},10^{-15}]$. We took $m$ as $10^6$, $n$ as $300$ and $k =2n = 600$. It turned out that the generated matrices $\bX^{(i)}$ had a rank of about $r  \approx 290$. The QR factorizations were calculated in float64 arithmetic. The $\tau$ parameter in RRRCholeskyQR was chosen to be $5\times 10^{-16}$. From~\cref{fig:Ex3_1a,fig:Ex3_1b} one can see that the standard Cholesky QR2 here failed completely, while the RCholeskyQR and the shifted CholeskyQR2 showed highly deteriorated stability than before. In contrast, the RRRCholeskyQR factorization still was almost perfectly stable, as was the Householder QR. This once again proves the unconditional stability of RRRCholeskyQR. It is important to note that it was revealed that the shifted CholeskyQR3 and RCholeskyQR2 not only could not provide an approximately orthonormal Q factor for large values of $\sigma$ (see~\cref{fig:Ex3_1c}), but they even failed several times due to the numerical indefiniteness of $\bQ^\mathrm{T} \bQ$ at the last CholeskyQR step, whereas RRRCholeskyQR2 showed great stability.
	
	\subsection{Randomized block GMRES}
	Next, we tested the stability of the column-oriented RCholeskyQR in the context of constructing a well-conditioned Krylov basis for solving a block linear system of equations
	$$\bA \bU = \bB,$$
	with randomized GMRES method~\cite{balabanov2021randomized,balabanov2021randomizedGS}.
	We have taken the linear system considered in numerical experiments in~\cite{balabanov2021randomized}. Namely, we took $\bA = (\bA_{Ga}+ 0.2 \bI)\bP_{Ga}$, where $\bA_{Ga}$ is the ``Ga41As41H72'' matrix of dimension $m = 268096$ from the { SuiteSparse matrix collection}, and $\bP_{Ga}$ is the incomplete LU preconditioner of $\bA_{Ga}+ 0.2 \bI$ with zero level of fill-in and symmetric reverse Cuthill-McKee reordering. {The matrix $\bA$ was not computed explicitly, but provided as an implicit map that outputs product with vectors and matrices}. The right hand side matrix $\bB$ was taken as an $m \times 100$ Gaussian matrix. This system has been approximately solved using the GMRES method based on various versions of the block Gram-Schmidt process or the column-oriented RCholeskyQR. We restarted GMRES every $30$ iterations, i.e. when the dimension of the Krylov space became $n = 3100$. In the randomized algorithms, the sketching dimension was chosen to be $k =7500$. In the experiments, the products with $\bA$ were calculated in float64 format. Solutions to Hessenberg least-squares problems in GMRES were obtained with Givens rotations that were also performed in float64 format. All operations related to the orthogonalization process were performed and accumulated in float32 format.
	
	In the column-oriented version of RCholeskyQR (\cref{alg:colRCholeskyQR}), we computed $\bR_{(1:i-1,i)}$ and $\bS_{(i)}, \bR_{(i , i) }$ in steps 3 and 4 with Householder QR.
	We here tested the variant of the block RGS algorithm that completely repeats the column-oriented RCholeskyQR (\cref{alg:colRCholeskyQR}) with updating the sketch $\bS_{(i-1)} \leftarrow \bTheta \bQ_{(i -1)}$ in step 2 as described in~\cref{colRCholeskyQR}. In the deterministic BCGS2 and BMGS algorithms, the inter-block orthogonalizations were performed by Householder QR.
	
\begin{figure}[!h]
		\centering
		\begin{subfigure}{.35\textwidth}
			\centering  
			\includegraphics[width=\textwidth]{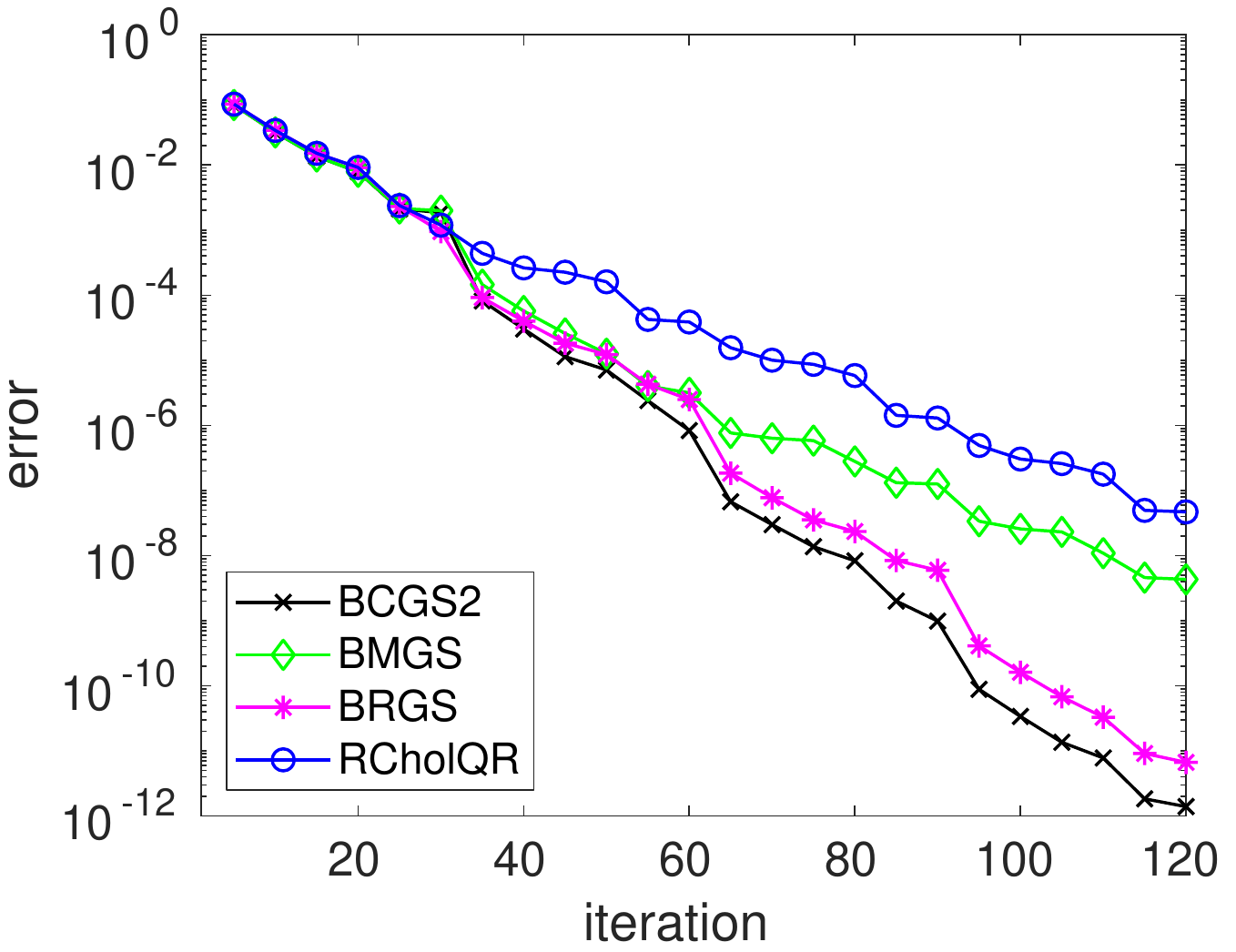}
			\caption{\small Max column-wise residual error \\
				\hfill  $~~~~~~~~~\max_{j}\|\bA \bU_{(j)} - \bB_{(j)} \|_2/\|\bB_{(j)}\|_2$.}
			\label{fig:Ex4_1a}
		\end{subfigure} \hspace{.03\textwidth}
		\begin{subfigure}{.35\textwidth}
			\centering
			\includegraphics[width=\textwidth]{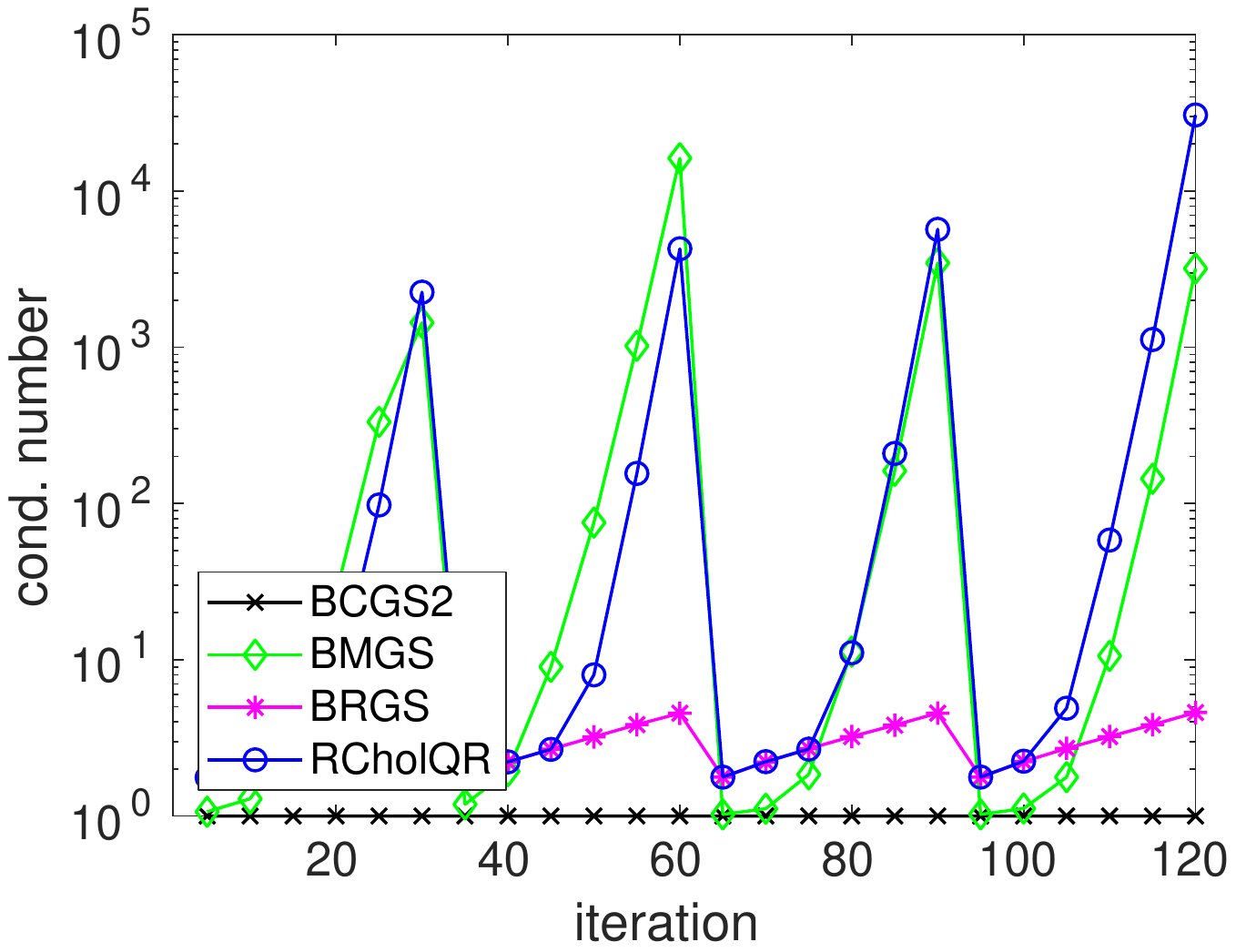}
			\caption{\small Cond. number of Krylov basis $\bQ_{(1:i)}$.}
			\label{fig:Ex4_1b}
		\end{subfigure}
		\caption{Solution of a linear system with GMRES.}
		\label{fig:Ex4_1}
	\end{figure}

	\Cref{fig:Ex4_1} depicts the convergence of the residual error $\max_{1\leq j \leq 100}\|\bA \bU_{(j)} - \bB_{(j)} \|_2/\| \bB_{(j)}\|_2$ and the condition number of the computed Krylov basis $\bQ_{(1:i)}$ at each iteration $i$. The BCGS2 remained perfectly stable throughout all iterations.  However, it is the most expensive algorithm tested, requiring nearly four times as many flops as the block RGS and RCholeskyQR algorithms. The block RGS also remained stable at all iterations. It provided a well-conditioned Q factor with $\mathrm{cond}(\bQ) \leq 5$ and almost as good residual error as BCGS2. We observe instabilities in BMGS and RCholeskyQR, which, unfortunately, worsened the convergence of the solution. In fact, it can be seen that RCholeskyQR entailed a residual error and the condition number of the Krylov basis that are almost four orders of magnitude greater than those of the block RGS. From this, we conclude that block RGS should be preferred over RCholeskyQR in the context of solving linear systems.

	\subsection{Runtime comparison}
	
	In this subsection, we explore the speedups that can be achieved with the proposed methods for computing QR factorizations with an orthonormal Q factor. The following experiments were performed in MATLAB R2021b on a node with 192GB of RAM and 2x Cacade Lake Intel Xeon 5218 16 cores 2.4GHz processor. We generated a sequence of random matrices $\bX = \bU \bSigma \bV^\mathrm{T}$ of different sizes and ranks, where $\bU$ and $\bV$ are orthonormalized random Gaussian matrices, $\bSigma = \mathrm{diag}([1,\sigma^{\frac {1}{r-1}}, \hdots, \sigma^{\frac{r-2}{r-1}},\sigma])$, and $\sigma=10^{-15}$. Then such $\bX$ were orthonormalized with RCholeskyQR2, RRRCholeskyQR2, shifted CholeskyQR3 and Householder QR. To compute the Householder QR we used the MATLAB's built-in function $\mathtt{qr}$. In shifted CholeskyQR3 we used the built-in $\mathtt{chol}$ function for the computation of Cholesky decomposition, and the built-in BLAS-3 forward substitution and matrix-matrix multiplication for other operations. In RCholeskyQR2 and RRRCholeskyQR2 algorithms we chose $\bTheta$ as a Gaussian OSE with twice as many rows as there are columns in $\bX$. 
	Furthermore, in randomized algorithms, in addition to the total runtime, we also measured the runtime corresponding to an ideal scenario where the sketching step requires negligible computational cost compared to other operations. In the given architecture, this could be achieved, for instance, with a well-implemented SRHT embedding.
	
\begin{figure}[!h] 
		\centering
		\begin{subfigure}{.35\textwidth}
			\centering  
			\includegraphics[width=\textwidth]{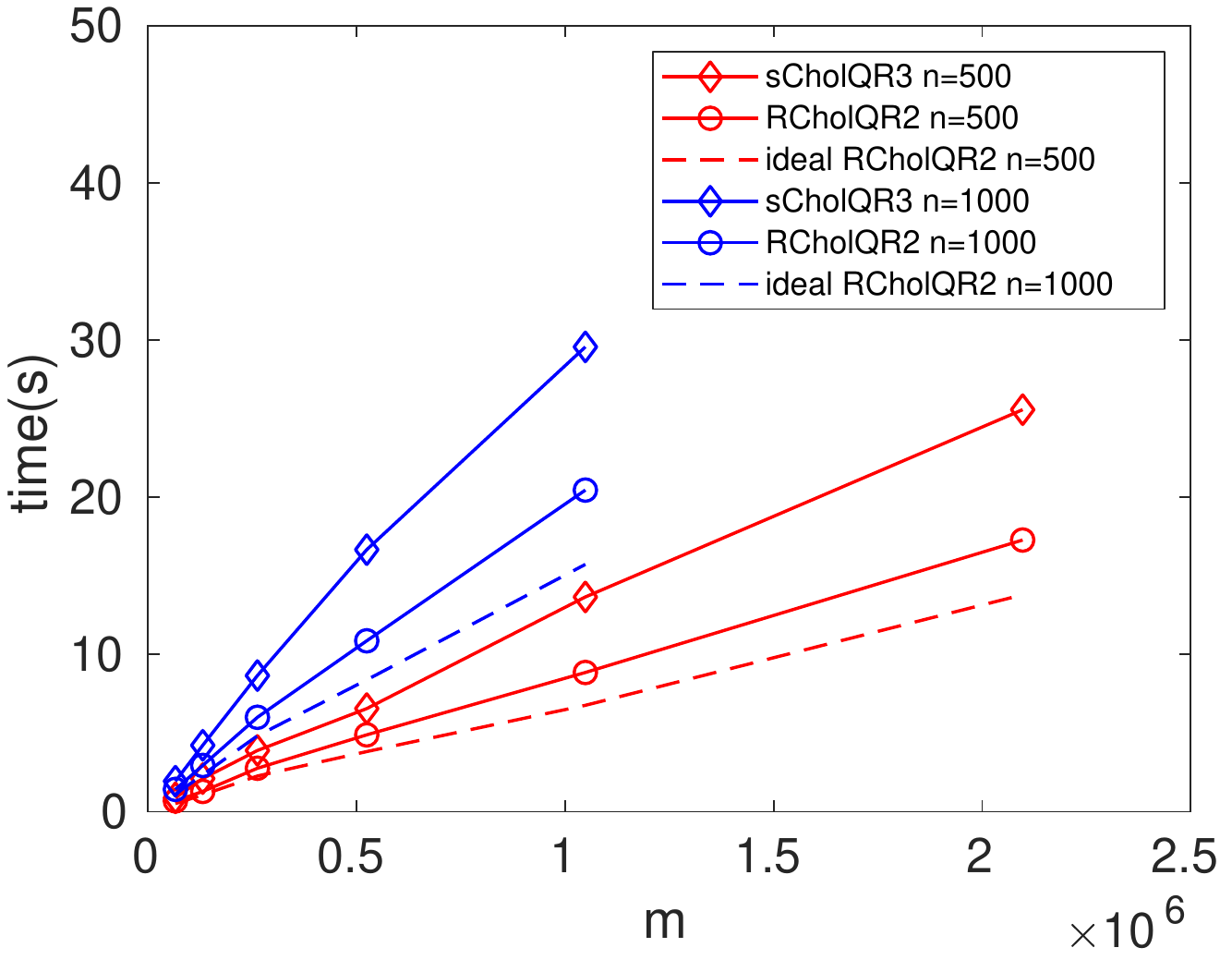}
			\caption{\small RCholQR2 vs. shifted CholQR3}
			\label{fig:Ex5_1a}
		\end{subfigure} \hspace{.03\textwidth}
		\begin{subfigure}{.35\textwidth}
			\centering
			\includegraphics[width=\textwidth]{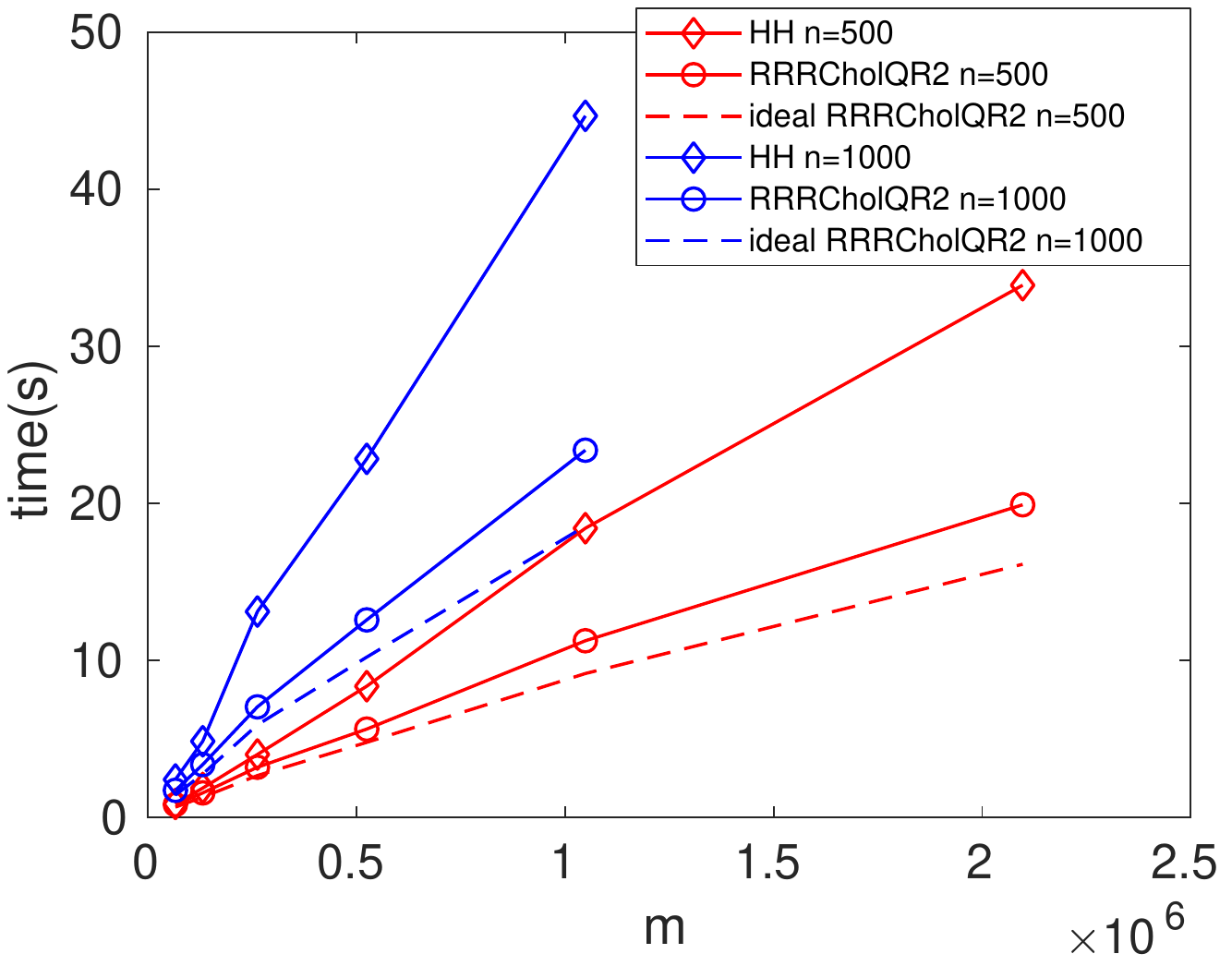}
			\caption{\small RRRCholQR2 vs. Householder}
			\label{fig:Ex5_1b}
		\end{subfigure}
		\caption{Runtimes in seconds taken by the QR factorizations of full-rank matrices $\bX$ of varying sizes.}
		\label{fig:Ex5_1}
	\end{figure}

	It can be seen from~\cref{fig:Ex5_1a} that the proposed RCholeskyQR2, despite being as or even more stable than the shifted CholeskyQR3, was up to $1.5$ times faster. In addition, this speedup could potentially be increased to a $2$ factor by using a more efficient sketching step. Let us now compare runtimes of the RRRCholeskyQR2 algorithm and Householder QR that are both unconditionally stable. We find from~\cref{fig:Ex5_1b} that for large $\bX^{(i)}$ RRRCholeskyQR2 required almost half as much runtime as Householder QR, which could potentially be reduced even more. 
	
	Furthermore, as was said, the RRRCholeskyQR2 algorithm not only provides the benefit of unconditional stability, but also has the ability to significantly reduce the computational cost when the matrix $\bX$ is of relatively low rank. This fact was validated too. From~\cref{tab:rankruntimes} we reveal that RRRCholeskyQR2 was almost $3.5$ times faster than the Householder QR when $\bX$ had a moderate rank, which could be improved to $10$ times, or potentially even $100$ times, when $\bX$ was of low rank. 
	
	\begin{table}[tbhp]
		\caption{Runtimes in seconds taken by QR factorizations of matrices $\bX$ of size $m=2^{20}$ and $n = 1000$, and of varying ranks $r$. For randomized algorithms, next to the overall runtimes we also provide the runtimes that could be achieved if the sketching step had a negligible computational cost. }
		\label{tab:rankruntimes}
		\centering
		\scalebox{0.9}{
			\begin{tabular}{|l|l|l|l|l|l|l|l|l|} 
				\hline
				& \multicolumn{2}{c|}{$r=1000$}  &\multicolumn{2}{c|}{$r=500$} & \multicolumn{2}{c|}{$r=100$} & \multicolumn{2}{c|}{$r=10$}  \\
				\hline
				{RCholQR2}  & $20.4$ & $15.7$ & $20.4$ & $15.7$  & $20.7$ & $16.1$ & $20.7$ & $16.1$ \\ [2pt] \hline
				{RRRCholQR2}  & $23.4$ & $18.5$ &  $14.1$  & $9.2$  & $6.2$ & $1.6$ & $4.4$ & $0.35$ \\ [2pt] \hline		
				{sCholQR3}   & $29.6$ &-- & $31.6$ & --  & $31.2$ & -- & $29.7$ & -- \\ [2pt] \hline		
				{HH}  & $44.6$ & -- &  $47.5$ & -- & $50.5$ & --&  $44.8$ & --   \\ [2pt] \hline
			\end{tabular}
		}
	\end{table}
	
	\section{Conclusion} \label{conclusion}
	
	This article proposed several variants of randomized Cholesky QR factorization. The presented direct RCholeskyQR algorithm should be up to four times more efficient than standard/shifted CholeskyQR2. Yet, it is just as or even more stable, and provides a well-conditioned Q factor whenever the input matrix is numerically full-rank. If necessary, RCholeskyQR can be readily augmented with the standard CholeskyQR to provide a Q factor that is orthonormal to machine precision, and not just well-conditioned, which results in RCholeskyQR2 algorithm. We have depicted some derivatives of RCholeskyQR, such as the column-oriented RCholeskyQR and the reduced RCholeskyQR. These algorithms can be useful for instance for constructing a Krylov basis, or for computing an approximation of a linear system's solution on a reduced basis. In addition, we have proposed an unconditionally stable RRRCholeskyQR, which can be seen as a very desirable alternative to other existing unconditionally stable algorithms such as Householder QR or TSQR. The RRRCholeskyQR should have the same computational cost as RCholeskyQR, or even better if the input matrix is of low rank.
	The proposed methodology was supported by rigorous theoretical and numerical stability analysis. The efficiency gains were also verified both theoretically and experimentally. In particular, in numerical experiments, we revealed speedups of RCholeskyQR2 and RRRCholeskyQR2 by almost $1.5$ and $2$ compared to the shifted CholeskyQR3 and Householder QR, respectively. These improvements could be made even greater by using a more efficient sketching subroutine. In addition, when the input matrix $\bX$ was of low rank, RRRCholeskyQR2 provided a much higher speedup, namely by a factor of $10$ or even $100$ (potentially).
	
	\section{Acknowledgments}
	The author would like to thank Laura Grigori for useful discussions on randomized algorithms.
	This project has received funding from the European Research Council (ERC) under the European Union's Horizon 2020 research and innovation program (grant agreement No 810367). 
	\printbibliography

	\end{document}